\documentclass[a4paper,10pt]{article}
\usepackage{amsfonts,amsmath,amssymb,enumitem}%,showkeys}

\newcommand{\CC}{\mathcal{C}}
\newcommand{\DD}{\mathcal{D}}
\newcommand{\EE}{\mathcal{E}}
\newcommand{\LL}{\mathcal{L}}
\newcommand{\UU}{\mathcal{U}}
\newcommand{\R}{\mathbb{R}}
\newcommand{\C}{\mathbb{C}}
\newcommand{\N}{\mathbb{N}}
\newcommand{\eps}{\varepsilon}
\newcommand{\om}{\omega}

\newcommand{\bsob}{\mathrm{H}_0}
\newcommand{\sob}{\mathrm{H}}
\newcommand{\Sa}{\mathrm{Sa}}
\newcommand{\spec}{\mathrm{spec}}
\newtheorem{Lemma}{Lemma}
\newtheorem{Proposition}{Proposition}
\newtheorem{Corollary}{Corollary}
\newtheorem{Definition}{Definition}
\newtheorem{Theorem}{Theorem}

\begin{document}

\title{A Convergence Result for Dirichlet Semigroups on Tubular Neighbourhoods and the Marginals of Conditional Brownian Motion}

\author{Vera Nobis \\ Lehrstuhl A f\"ur Mathematik \\ Olaf Wittich \\ Lehrstuhl I f\"ur Mathematik \\
RWTH Aachen University}

\maketitle

\abstract{We investigate yet another approach to understand the limit behaviour of Brownian motion conditioned to stay within a tubular neighbourhood around a closed and connected submanifold of a Riemannian manifold. In this context, we identify a second order generator subject to Dirichlet conditions on the boundary of the tube and study its associated semigroups. After a suitable rescaling and renormalization procedure, we obtain convergence of these semigroups, both in $L^2$ and in Sobolev spaces of arbitrarily large index, to a limit semigroup, as the tube diameter tends to zero. As a byproduct, we conclude that the conditional Brownian motion converges in finite dimensional distributions to a limit process supported by the path space of the submanifold. \\


}
	
\section{Introduction}\label{intro}

\noindent We consider Brownian motion on a complete Riemannian manifold $M$ conditioned not to leave a tube $L(\eps)$ of small radius $\eps >0$ around a closed and connected submanifold $L$. We ask the question, whether a sequence of path measures obtained in this way, converges weakly to a measure supported by the path space of the submanifold as $\eps$ tends to zero. This question was answered to the affirmative for embeddings into Euclidean space in \cite{SidSmoWeiWit:03a} using methods from stochastic differential equations. For embeddings into general Riemannian manifolds, we follow a different approach via a perturbational ansatz. Starting from the connection between conditioned and absorbed process explained in \ref{proc} below, we identify a second order generator $H_{\eps}$ subject to Dirichlet boundary conditions. The associated semigroups are transformed to a fixed tube $L(1)$ and suitably renormalized. They correspond to the one-dimensional marginals of conditional Brownian motion transformed by a multiplicative functional. The convergence result Theorem \ref{Main} for the semigroups implies convergence of the associated processes in finite dimensional distributions. In a subsequent paper \cite{Wit:19}, we will prove that this sequence of measures is actually tight, which even implies weak convergence of the path measures.\\  

\noindent The paper is organized as follows: First we introduce the main result Theorem \ref{Main} below, and explain why it implies convergence of the associated processes in finite - dimensional distributions. In Section \ref{perturbation}, we give a precise description of the perturbation problem under consideration and of the Sasaki metric on the tube. Assuming some knowledge about the terms in the perturbation expansion from Proposition \ref{pertprop}, we prove epi-convergence of the quadratic forms associated to the generators and conclude convergence of the semigroups in an $L^2$-sense. In Section \ref{geometry}, we investigate the geometry of small tubes around submanifolds. In particular, we compare the induced metric with the Sasaki metric on the tube and prove Proposition \ref{pertprop}. Since $L\subset M$ is a zero set, $L^2$-convergence of the semigroups is not sufficient to prove convergence of the conditional process. Therefore, in the final section, we establish some a priori estimates for analytic vectors in the domain of the generators and use them to finally show that the semigroups actually converge smoothly in the sense of Theorem \ref{Main} below. \\

\noindent Please note that, if not indicated otherwise, $\Vert - \Vert$ denotes the norm on $L^2 (L(1),\mu_{\Sa})$ and $\Vert - \Vert_n$ the norm on the Sobolev space $\sob^n (L(1),\mu_{\Sa})$.

\subsection{A family of semigroups and its convergence}

\noindent a. Let $L\subset M$ be a closed Riemannian submanifold of the Riemannian manifold $M$. We assume, without loss of generality, that the exponential map maps a neighbourhood of the zero section in the normal bundle $NL$ diffeomorphically onto the $r$-tube $L(r):= \lbrace x\in M\,:\,d_M(x,L)<r\rbrace$ for some $r>1$. On $L(r)$, we study two different metrics, the metric $g$ induced by the embedding into $M$ and the {\em Sasaki metric} $g_{\Sa}$. Since $L(1)$ with either metric is assumed to be diffeomorphic to the unit disc bundle $D_1NL :=\lbrace W\in NL\,:\,\Vert W\Vert < 1\rbrace$, of the normal bundle, we will denote both spaces by $L(1)$ without further mentioning. Moreover, we will not distinguish between the respective metrics on $L(1)$ and its pullbacks to the disc bundle. Let $\star$ denote the {\em Hodge operator} associated to $g$, and $\rho := \frac{d\mu}{d\mu_{\Sa}}$ the Radon - Nikodym density of the volume forms $\mu$ and $\mu_{\Sa}$ of the respective metrics. The Laplace-Beltrami operator on $M$ is denoted by $\Delta:= -d^{\star}d\geq 0$, the tube projection by $\pi:L(1)\to L$ and, for $\eps > 0$, the {\em rescaling map} $\sigma_{\eps} : L(\eps)\to L(1)$ is given by $\sigma_{\eps}(u):= \exp_{\pi u}(\eps^{-1}\exp^{-1}_{\pi u}(u))$. \\

\noindent b. Let now $0<\eps \leq 1$ and consider the smooth potential 
\begin{equation}\label{potential}
U:=\rho^{1/2}\Delta\rho^{-1/2}\in C^{\infty}(L(r)).
\end{equation} 
By $(\Delta-U)_{\eps}$, we denote the Hamiltonian on $L(\eps)$ with Dirichlet boundary conditions on $\partial L(\eps)$, i.e. the operator associated to the quadratic form
$$
Q_{\eps}:\bsob^1(L(\eps),\mu)\to\R, \,\, f\mapsto \int_{L(\eps)} df\wedge \star df - \star\,Uf^2
$$
by Friedrichs' construction.  \\

\noindent c. Let now $\Sigma_{\eps}:L^2(L(1),\mu_{\Sa})\to L^2(L(\eps),\mu)$ be the map given by
\begin{equation}\label{rescaling_map}
\Sigma_{\eps}f:=\left(\eps^{m-l}\rho\right)^{-1/2} \sigma_{\eps}^*f.
\end{equation}
By partial integration, it turns out that 
\begin{equation}\label{laplacemainops}
\Sigma_{\eps}^{-1}\circ(\Delta-U)_{\eps}\circ\Sigma_{\eps} = H_{\eps},
\end{equation}
where $H_{\eps}$ with domain $\bsob^1\cap\sob^2(L(1),\mu_{\Sa})$ is self-adjoint and non-negative on $L^2(L(1),\mu_{\Sa})$.\\

\noindent Because the parameter $\eps > 0$ is closely related to the tube radius, the perturbation problem for $H_{\eps}$ is not to be expected to yield a sensible limit as $\eps$ tends to zero. However, if $\lambda_0>0$ is the smallest eigenvalue for the Dirichlet problem on the $(m-l)$-dimensional Euclidean unit ball $B\subset\R^{m-l}$, the semigroups generated by $H_{\eps}^0 := H_{\eps} - \eps^{-2}\lambda_0$ will converge strongly to a semigroup on a certain subspace $E_0\subset L^2(L(1),\mu_{\Sa})$ (cf. \ref{sasaki_can_var}.b). Denoting the orthogonal projection onto the subspace $E_0$ by the same symbol, the main result of this paper reads as follows:

\begin{Theorem}\label{Main} Let $u(\eps)\in L^2(L(1),\mu_{\Sa})$ be a strongly continuous family of functions and denote by $\Delta_L$ the Laplace-Beltrami operator on $L$. Then, for all $n\geq 1$, we have
	\begin{equation}\label{result}
	\lim_{\eps\to 0} e^{-\frac{t}{2}H_\eps^0}u(\eps) = E_0\,e^{-\frac{t}{2}\Delta_L} \,E_0u(0) 
	\end{equation}
	uniformly on each compact sub-interval $I\subset (0,\infty)$ in the Sobolev space $\sob^{2n}(L(1),\mu_{\Sa})$.  
\end{Theorem}

\noindent{\bf Remark.} (a) From Proposition \ref{pertprop} (1) below, we obtain for an arbitrary $f\in C^{\infty}(M)$
$$
E_0 f = \phi_0 \,\left(\langle f,\phi_0\rangle_{\pi^{-1}(x)}\right)\circ\pi,
$$
where $\langle -,-\rangle_{\pi^{-1}(x)}$ is the scalar product with respect to the Riemannian volume on the fibre induced by the Sasaki metric and $\phi_0\in L^2(L(1),\mu_{\Sa})$ is an explicitly given function. (For the precise definition of $\phi_0$ see Section \ref{proof_pertprop}.) The precise meaning of the right hand side in (\ref{result}) is therefore given by
\begin{equation}\label{meaning}
E_0\,e^{-\frac{t}{2}\Delta_L} \,E_0f = \phi_0\,\left(e^{-\frac{t}{2}\Delta_L}f_a\right)\circ\pi,
\end{equation}
where $f_a\in C^{\infty}(L)$ is given by $x\mapsto \langle \phi_0,f\rangle_{\pi^{-1}(x)}$. (b) Theorem \ref{Main} will still hold if $u(\eps)$ is only strongly continuous at $\eps=0$.

\subsection{A corresponding conditional process and its convergence}\label{proc}

\noindent Theorem \ref{Main} is related to the fact that Brownian motion on $M$ conditioned to smaller and smaller tubes $L(\eps)$ around $L$, converges to a process with a path measure which is equivalent to the Wiener measure on $L$. For embeddings into Euclidean space this was shown in \cite{SidSmoWeiWit:03b}. In this section, we are going to discuss this connection.\\

\noindent a. Let $\Omega := C(\lbrack 0,\infty), M)$ be the path space and 
$$
\Omega^{\eps}_{s,t}:= \lbrace \om\in\Omega\,:\,\om(u)\in L(\eps), s\leq u \leq t\rbrace.
$$  
Denoting by $\mathbb{W}$ the Wiener measure on $M$, we fix some finite $T>0$ and consider the measure 
$$
\nu (d\om) = \exp\left(\int_0^{T}\varphi U(\om(s)) \,ds\right)\,\mathbb{W}(d\om),
$$
on the path space of $M$. Here, $\varphi\in C^{\infty}(M)$ is smooth, $\varphi\vert_{L(1)} = 1$, and $\varphi\vert_{M\setminus L(r)}=0$. Now we consider the probability measures $\mathcal{L}_{\eps}$, $\eps > 0$, which are obtained by restricting $\nu$ to the set $\Omega^{\eps}_{0,T}$ followed by normalization to total mass one. To be precise, 
\begin{equation*}
\mathcal{L}_{\eps} (d\om) = \frac{\nu (d\om\cap \Omega^{\eps}_{0,T})}{\nu(\Omega^{\eps}_{0,T})}
\end{equation*}
supported by the path space of $L(\eps)$. The processes with distribution $\mathcal{L}_{\eps}$ are denoted by $(x_t^{\eps})_{0\leq t\leq T}$. For $0\leq s < t\leq T$, we define the transition kernel of $x_t^{\eps}$ by the conditional probability
\begin{equation*}
 Q_{\eps}(s,x;t,dy) := \mathcal{L}_{\eps} (\om (t)\in dy\,\vert\, \om (s) = x).
\end{equation*}

\noindent b. By the Markov property of Wiener measure and the properties of conditional expectation that implies
\begin{eqnarray*}
Q_{\eps}(s,x;t,dy) &=& \nu (\om (t)\in dy\,\vert\, \om\in\Omega^{\eps}_{0,T},\om(s) = x)\\
&=& \frac{\nu(\om (t)\in dy, \om\in \Omega^{\eps}_{s,T}\,\vert\,\om (s) = x)}{\nu (\Omega^{\eps}_{s,T}\,\vert\, \om(s) = x)}\\
&=& \frac{\nu (\Omega_{t,T}^{\eps}\,\vert\, \om(t) = y)\nu(\om(t)\in dy, \om\in\Omega_{s,t}^{\eps}\,\vert\,\om (s) = x)}{\nu (\Omega^{\eps}_{s,T}\,\vert\, \om(s) = x)}.
\end{eqnarray*}
The crucial observation which establishes the connection between the conditional process and the Dirichlet operator considered above is now that 
\begin{eqnarray*}
P_{\eps}(s,x;t,dy) &=& \nu (\om(t)\in dy, \om\in\Omega_{s,t}^{\eps}\,\vert\,\om (s) = x) \\
&=& \nu(\om(t)\in dy, t < \tau_{\eps}(\om)\,\vert\,\om (s) = x),
\end{eqnarray*}
where $\tau_{\eps}$ is the first exit time from $L(\eps)$ and hence,
\begin{equation}\label{conditional_kernel}
 Q_{\eps}(s,x;t,dy) = \frac{\pi_{T-t}^{\eps}(y)}{\pi_{T-s}^{\eps}(x)}P^{\eps}(s,x;t,dy)
\end{equation}
where $\pi^{\eps}_{u}(w) :=\int_{L(\eps)}P^{\eps}(0,w;u,dz)$.\\

\noindent c. By the Feynman-Kac formula, integration with respect to the transition kernel can be represented probabilistically by
\begin{eqnarray*}
& & \int_{L(\eps)} f(y)P^{\eps}(s,x;t,dy) \\ &=& \int_{\Omega} f(\omega(t))\nu(\om(t)\in dy, t < \tau_{\eps}(\om)\,\vert\,\om (s) = x) \\
&=& \int_{\Omega} f(\omega(t))\exp\left(\frac{1}{2}\int_s^{t} U(\om(s)) ds\right)\mathbb{W}(\om(t)\in dy, t < \tau_{\eps}(\om)\vert\om (s) = x), 
\end{eqnarray*}
and in terms of generators and semigroups, we have
$$
\int_{L(\eps)} f(y)P^{\eps}(s,x;t,dy) = \left(e^{-\frac{t-s}{2}(\Delta-U)_{\eps}}f\right)(x).
$$
Hence, we obtain from (\ref{conditional_kernel}) for the conditional process starting at $s=0$ in $x\in L$ that
\begin{equation}\label{semigrouprep}
E^x\lbrack f(x_t^{\eps})\rbrack = \frac{e^{-\frac{t}{2}(\Delta-U)_{\eps}}\left(fe^{-\frac{T-t}{2}(\Delta-U)_{\eps}}1\right)}{e^{-\frac{T}{2}(\Delta-U)_{\eps}}1}(x).
\end{equation}

\noindent d. In this subsection, we are going to explain how we can use Theorem \ref{Main} together with the statements (\ref{laplacemainops}) and (\ref{semigrouprep}), to conclude that the processes $x_t^{\eps}$ converge to Brownian motion $x^0_t$ on $L$ in finite dimensional distributions. For Markov processes, the following statement about convergence of the one-dimensional marginals implies convergence in finite dimensional distributions.

\begin{Corollary}\label{findim_conv} Let $x^0_{t, 0\leq t\leq T}$ be Brownian motion on $L$. Let $x_0^{\eps} = x\in L$ be a fixed common starting point. Then, for all $f\in C^{\infty}(M)$ and $0\leq t\leq T$, we have 
$$
\lim_{\eps\to 0} 
E^x\lbrack f(x_t^{\eps})\rbrack = E^x\lbrack f\vert_L(x_t^0)\rbrack,
$$
i.e. the associated flows converge as $\eps$ tends to zero. 
\end{Corollary}

\noindent{\bf Proof} From (\ref{semigrouprep}), using the rescaling map $\Sigma_{\eps}$ from \ref{intro}.c together with $\sigma_{\eps}^*(f) = f\circ\sigma_{\eps}$, $\sigma_{\eps}(x) = x$, we obtain 
\begin{eqnarray*}
E^x\lbrack f(x_t^{\eps})\rbrack &=& \frac{\Sigma_{\eps}e^{-\frac{t}{2} H_{\eps} }\Sigma_{\eps}^{-1}\left(f\Sigma_{\eps} e^{-\frac{T-t}{2}H_{\eps}}\Sigma_{\eps}^{-1}1\right)}{\Sigma_{\eps} e^{-\frac{T}{2}H_{\eps}}\Sigma_{\eps}^{-1}1}(x)\\
&=&\frac{e^{-\frac{t}{2}H_{\eps}^0}\left(\sigma_{\eps}^{*\,-1}(f)e^{-\frac{T-t}{2}H_{\eps}^0}\sigma_{\eps}^{*\,-1}(\sqrt{\rho})\right)}{e^{-\frac{T}{2}H_{\eps}^0}\sigma_{\eps}^{*\,-1}(\sqrt{\rho})}(x).
\end{eqnarray*}
Now, $\sigma_{\eps}^{*\,-1}(\sqrt{\rho})\to 1$ uniformly on $L(1)$ as $\eps$ tends to zero, since $\rho\in C^{\infty}(L(1))$ is smooth and $\rho\vert_L = 1$. On the other hand $\sigma_{\eps}^{*\,-1}(f)\to f\vert_L\circ\pi$. This fact combined with Theorem \ref{Main} yields
$$
u(\eps) = \sigma_{\eps}^{*\,-1}(f)\,e^{-\frac{T-t}{2}H_{\eps}^0}\sigma_{\eps}^{*\,-1}(\sqrt{\rho})
$$
is strongly continuous, and therefore, again by Theorem \ref{Main}, the right hand side above converges to
$$
P_tf :=\frac{E_0e^{-\frac{t}{2}\Delta_L}E_0\left((f\vert_L\circ\pi) \, E_0e^{-\frac{T-t}{2}\Delta_L}E_01\right)}{E_0e^{-\frac{T}{2} \Delta_L}E_01}.
$$ 
By (\ref{meaning}) and by $\pi (x) = x$, that finally implies 
\begin{equation*}
P_tf = \left(\frac{e^{-\frac{t}{2}\Delta_L}\left((f\vert_L) e^{-\frac{T-t}{2}\Delta_L}1\right)}{e^{-\frac{T}{2} \Delta_L}1}\right)\circ\pi(x) = E^x\lbrack f\vert_L(x_t^0)\rbrack.
\end{equation*}
\hfill $\Box$

\noindent Convergence of the marginals is the first part of proving weak convergence of the path measures. The second part is tightness of the measure family. Tightness of the measure family will be discussed in a subsequent paper.

\section{The Perturbation Problem}\label{perturbation}

\noindent First of all, we give a precise description of the Sasaki metric and an alternative description of the quadratic form associated to the operator $H_{\eps}$.\\

\noindent The Sasaki metric (cf. \cite{Reck:79}) on the normal bundle is given by  
$\langle X,Y\rangle_{\Sa} = \langle \pi_*X,\pi_*Y\rangle_L + \langle KX,KY\rangle_{NL}$ where $K$ denotes the connection map of the induced connection on the normal bundle and $\langle-,-\rangle_L$, $\langle-,-\rangle_{NL}$ the scalar product on $L$ and $NL$, respectively. For the cotangent bundle that implies for $W\in NL$
$$
\langle \xi,\eta\rangle_{\Sa, W} = \langle J_W^* \xi, J_W^* \eta\rangle_{N^*L} +\langle \kappa_W \xi, \kappa_W \eta\rangle_{L},
$$
where $J_W:N_{\pi W}L\to T_WN_{\pi W}L$ denotes the canonical isomorphism, $J_W^*$ its dual and $\kappa_W : T^*_WNL\to T^*_{\pi(W)}L$ is the dual of the horizontal lift for vector fields. \\

\noindent With these notations, the operator $\#_{\eps}: \wedge^{1}T^*L(1)\to \wedge^{m-1}T^*L(1)$ is given by
\begin{equation}\label{scalestar}
\#_{\eps} := \eps^{l-m}\left(\sigma_{\eps}^*\right)^{-1}\circ (\rho^{-1}\star)\circ \sigma_{\eps}^* ,
\end{equation}
with an associated bilinear form 
\begin{equation}\label{scaleform}
b_{\eps} (h,f) = \int_{L(1)} dh \wedge \#_{\eps} df .
\end{equation}
On $L^2 (L(1),\mu_{\Sa})$, the quadratic form
\begin{equation}\label{scalequad}
q_{\eps} (f) = \int_{L(1)} df \wedge \#_{\eps} df
\end{equation}
with domain $\DD = \bsob^1 (L(1),\mu_{\Sa})$ for all $\eps > 0$ is closed, non-negative and densely defined. By Friedrichs' construction, there is a self-adjoint and non-negative operator $H_{\eps}$ on $L^2(L(1),\mu_{\Sa})$ with domain $\DD (H_{\eps})=\bsob^1\cap\sob^2 (L(1),\mu_{\Sa})$ such that
\begin{equation}
b_{\eps} (h,f) = \int_{L(1)} h\,H_{\eps}f \,d\mu_{\Sa}.
\end{equation}
By $d(h\#_{\eps}df) = dh\wedge \#_{\eps}df + h d\#_{\eps}df$ and Stokes' theorem, the differential expression for $H_{\eps}$ is given by
\begin{equation}
H_{\eps}f = - \star_{\Sa} d\#_{\eps}df,
\end{equation}
where $\star_{\Sa}$ denotes the Hodge operator associated to $g_{\Sa}$. If the metric $g$ equals the Sasaki metric $g_{\Sa}$ (meaning in particular $\rho = 1$), the associated forms and operators will be denoted by $\#_{\Sa,\eps}$, $b_{\Sa,\eps}$, $q_{\Sa,\eps}$ and $H_{\Sa,\eps}$, respectively. \\

\noindent In the sequel, the quadratic form associated to the induced metric will be considered as a perturbation of the one associated to the Sasaki metric.

\subsection{Sasaki Metric and Canonical Variation}\label{sasaki_can_var}

\noindent a. The quadratic form $q_{\Sa,\eps}$ is nothing but the quadratic form of the Laplace - Beltrami operator associated to the {\em canonical variation} (cf. \cite{BerBou:82}, (5.1), p. 191) of the Sasaki metric. By
\begin{equation}
\eps^{m-l}\, d\mu_{\Sa} = d\mu_{\Sa}\circ\sigma_{\eps}^{-1}, \quad \lbrack J,\sigma_{\eps}^*\rbrack = \eps^{-1}\,\,\mathrm{and}\,\,\quad \lbrack \kappa,\sigma_{\eps}^*\rbrack = 0,
\end{equation}
we obtain
\begin{eqnarray*}
q_{\Sa,\eps} (f) &=& \int_{L(1)} \langle\sigma_{\eps}^*df, \sigma_{\eps}^*df\rangle_{\Sa}\circ\sigma_{\eps}^{-1} \,d\mu_{\Sa} \\
&=& \int_{L(1)} \left(\langle J^* \sigma_{\eps}^*df, J^* \sigma_{\eps}^*df\rangle_{N^*L} +\langle \kappa \sigma_{\eps}^*df, \kappa \sigma_{\eps}^*df\rangle_{L}\right)\circ\sigma_{\eps}^{-1} \,d\mu_{\Sa}\\
&=&\int_{L(1)} \left(\frac{1}{\eps^2}\langle J^* df, J^* df\rangle_{N^*L} +\langle \kappa df, \kappa df\rangle_{L}\right) \,d\mu_{\Sa},
\end{eqnarray*}
and the form can be written as a sum of two densely defined, closed quadratic forms on $L^2(L(1),\mu_{\Sa})$. 

\begin{Definition}\label{qvqh} {\rm (i)} The {\em vertical form} is given by
\begin{equation*}
q_{V}(f) := \int_{L(1)}\langle J^* df, J^* df\rangle_{N^*L}d\mu_{\Sa}, 
\end{equation*}
with domain
$$
\mathcal{D}_V :=\lbrace f\in L^2 (L(1),\mu_{\Sa})\,:\, f\vert_{\pi^{-1}(x)}\in\bsob^1(\pi^{-1}(x)) \,\mu_L-\mathrm{a.e.}, q_V(f) <\infty\rbrace,
$$
where $\mu_L$ denotes the Riemannian volume measure of the submanifold and the fibres are equipped with the metric induced by the Sasaki metric. {\rm (ii)} The {\em horizontal form} is given by
\begin{equation*}
q_{H}(f) := \int_{L(1)} \langle \kappa df, \kappa df\rangle_{L}d\mu_{\Sa}, 
\end{equation*}
with domain
$$
\mathcal{D}_H :=\lbrace f\in L^2 (L(1),\mu_{\Sa})\,:\, f\vert_{\partial L(r)}\in\mathcal{D}_r \,\,\lambda-\mathrm{a.e.}, q_H(f) <\infty\rbrace,
$$
where $\lambda$ denotes the Lebesgue measure on $(0,1)$ and the tube boundary $\partial L(r)$ is again equipped with the metric induced by the Sasaki metric and the induced Riemannian measure $\mu_r$. For $r\in (0,1)$, the domain $\mathcal{D}_r$ is given by
$$
\mathcal{D}_r:=\lbrace f\in L^2 (\partial L(r),\mu_r)\,:\, Xf\in L^2 (\partial L(r),\mu_r) \,\mathrm{for\, all} \, X \in \mathrm{Hor}_r \rbrace,
$$
where $\mathrm{Hor}_r$ denotes the space of smooth vector fields on $\partial L(r)$ with values in the restriction $H\vert_{\partial L(r)}$ of the horizontal subbundle $H\subset TNL$.
\end{Definition}
By $\DD = \DD_V\cap\DD_H$ we have indeed $q_{\Sa,\eps} = \eps^{-2} q_V + q_H$. The self - adjoint differential operators associated to the respective  quadratic forms by Friedrichs' construction are called {\em vertical and horizontal Laplacian} (cf. \cite{BerBou:82}, (1.2), p. 183). We denote the vertical operator by $\Delta_V$ and the horizontal operator by $\Delta_H$. Details of this construction are provided in Section \ref{forms_and_ops}.\\

\noindent b. If the total space $L(1)$ is equipped with the Sasaki metric, the projection $\pi:L(1)\to L$ will be a Riemannian submersion with totally geodesic fibres. Therefore, the fibres are isometric (\cite{Her:60}, 4.1). The prototype is the flat unit disc $B\subset\R^{m-l}$. Hence, all Dirichlet Laplacians 
$$
\Delta_x : \bsob^1\cap\sob^2 (\pi^{-1}(x))\to L^2(\pi^{-1}(x))
$$ 
on the fibres $\pi^{-1}(x)$ are unitarily equivalent with eigenvalues $0< \lambda_0 < \lambda_1 < ...$ and corresponding eigenprojections $E_{k,x}$. Therefore, applying Friedrichs' construction fibrewise, the vertical operator is given by a constant fibre direct integral (\cite{ReeSim:78}, p. 283) 
\begin{equation*}
\Delta_V = \int^{\oplus} \Delta_x,
\end{equation*}
with domain
$$
\DD (\Delta_V) := \left\lbrace f\in \int_L^{\oplus}\bsob^1\cap\sob^2(\pi^{-1}(x))d\mu_L\,:\, \, \int_L\,\Vert \Delta_x f\vert_{\pi(x)} \Vert_x^2 d\mu_L<\infty\right\rbrace.
$$
Hence, by \cite{ReeSim:78}, Theorem XIII.85, $\Delta_V$ is self - adjoint with a spectrum consisting precisely of the same eigenvalues $0< \lambda_0 < \lambda_1 < ... $ with corresponding eigenprojections
\begin{equation}
E_k := \int^{\oplus} E_{k,x} .
\end{equation}
In the sequel, we denote projections and corresponding eigenspaces by the same symbol.\\

\noindent c. The operator associated to $q_{\Sa,1}$ is the Dirichlet Laplacian $\Delta_{\Sa}$ on $L(1)$ with the Sasaki metric $g_{\Sa}$. $\Delta_{\Sa}$ is therefore the operator of a regular elliptic boundary value problem. By the compactness of $\overline{L(1)}$ (\cite{Tay:99}, 5.1, p. 303 ff.), the spectrum is discrete and consists only of eigenvalues $0 < \mu_{0} < \mu_{1} < ...$ of finite multiplicity. Furthermore, all eigenfunctions $\phi\in C^{\infty} (L(1))\cap C(\overline{L(1)})$ are continuous, smooth in the interior of the tube and vanish on the boundary. By \cite{BerBou:82}, (1.5), the operators $\Delta_{\Sa}$, $\Delta_V$, and $\Delta_H$ commute pairwise, meaning that in particular
\begin{equation*}
\lbrack \Delta_{\Sa},\Delta_V\rbrack f = \lbrack \Delta_{\Sa},\Delta_H\rbrack f = \lbrack\Delta_V,\Delta_H\rbrack f = 0
\end{equation*}
for all $f\in C^4(L(1))$. Thus, diagonalizing the operators $\Delta_H$ and $\Delta_V$ simultaneously on the finite dimensional eigenspaces of $\Delta_{\Sa}$, we obtain a common orthonormal base $u_{l}$, $l\geq 0$, of $L^2(L(1),\mu_{\Sa})$ of smooth eigenfunctions with
\begin{equation}
\Delta_{\Sa}u_{l} = \mu_{l}u_{l}, \, \Delta_Vu_{l} = \lambda_{k(l)}u_{l},  \,\Delta_Hu_{l} = (\mu_{l}- \lambda_{k(l)})u_{l},
\end{equation}
where $k(l)\in \N_0$ is the unique number such that $\Delta_Vu_{l} = \lambda_{k(l)}u_l$. Therefore, the operators commute as self - adjoint operators on $L^2(L(1),\mu_{\Sa})$. By the decomposition of $q_{\Sa,\eps}$ above, we obtain $\Delta_{\Sa} = H_{\Sa,1} = \Delta_V + \Delta_H$ and the spectral decomposition
\begin{equation}\label{saspec}
H_{\Sa,\eps} = \frac{1}{\eps^2} \Delta_V + \Delta_H = \sum_{l\geq 0} \left(\frac{\lambda_{k(l)}}{\eps^2} + \mu_{l}- \lambda_{k(l)}\right)u_l u_l^+.
\end{equation}
In particular, all Laplacians $H_{\Sa,\eps}$ associated to the canonical variation share the same eigenfunctions.\\

\noindent d. The operators $H_{\Sa,\eps}$, $\Delta_V$, and $\Delta_H$ are non-negative and self-adjoint. Therefore, they generate strongly continuous semigroups of contractions. By (c.) above, these semigroups commute pairwise. A suitably renormalized version of the semigroup generated by $H_{\Sa,\eps}$ converges strongly as $\eps$ tends to zero.

\begin{Lemma}\label{sasaconv} As $\eps$ tends to zero, we have 
\begin{equation*}
\lim_{\eps\to 0}e^{-\frac{t}{2}\left(H_{\Sa,\eps}-\frac{\lambda_0}{\eps^2}\right)}f = E_0 e^{-\frac{t}{2}\Delta_H}E_0f .
\end{equation*}
for all $f\in L^2(L(1),\mu_{\Sa})$.
\end{Lemma}

\noindent{\bf Proof.} The semigroups generated by $\Delta_V$ and $\Delta_H$ commute, hence by the spectral theorem
\begin{eqnarray*}
 e^{-\frac{t}{2}\left(H_{\Sa,\eps}-\frac{\lambda_0}{\eps^2}\right)}f &=& e^{-\frac{t}{2}\left(\frac{\Delta_V-\lambda_0}{\eps^2}\right)}e^{-\frac{t}{2}\Delta_H}f \\
 &=& \sum_{k\geq 0} e^{-\frac{t(\lambda_{k}-\lambda_0)}{2\eps^2}}E_ke^{-\frac{t}{2}\Delta_H}f \\
 &=& E_0e^{-\frac{t}{2}\Delta_H}E_0f + \sum_{k\geq 1} e^{-\frac{t(\lambda_{k}-\lambda_0)}{2\eps^2}}E_ke^{-\frac{t}{2}\Delta_H}f.
\end{eqnarray*}
By contractivity of the semigroup, we have
$$
\left\Vert\sum_{k\geq 1} e^{-\frac{t(\lambda_{k}-\lambda_0)}{2\eps^2}}E_ke^{-\frac{t}{2}\Delta_H}f\right\Vert \leq e^{-\frac{t(\lambda_{1}-\lambda_0)}{2\eps^2}}\,\Vert e^{-\frac{t}{2}\Delta_H}f\Vert\leq e^{-\frac{t(\lambda_{1}-\lambda_0)}{2\eps^2}}\,\Vert f\Vert
$$
and this tends to zero as $\eps$ tends to zero.\hfill $\Box$ \\

\noindent{\bf Remark.} The action of the semigroup generated by $\Delta_H$ on $E_0$ will be described more explicitly in Section \ref{geometry}.

\begin{Definition} In the sequel, the objects
\begin{equation*}
\begin{array}{ll}
H_{\eps}^0 := H_{\eps} - \eps^{-2}\lambda_0, & q_{\eps}^0(f) := q_{\eps}(f) - \eps^{-2}\lambda_0\Vert f\Vert^2
\end{array}
\end{equation*}
are called {\em renormalized operator}, and {\em renormalized form}, respectively. 
\end{Definition}
The following inequality will be very helpful to understand the perturbation and follows from the spectral properties considered above.

\begin{Lemma}\label{est1} (i) There are constants $a,A > 0$ such that
\begin{equation*}
a\, q_{\Sa,1}(f) \leq \Vert f \Vert_{\bsob^1(L(1),\mu_{\Sa})}^2 \leq A \, q_{\Sa,1}(f) 
\end{equation*}
for all $f\in \bsob^1 (L(1),\mu_{\Sa})$. (ii) There is a constant $k_{\Sa} \geq 1$ such that
\begin{equation*}
q_{V}(f) \leq k_{\Sa} \left( \eps q_{\Sa,\eps}^0(f) + \Vert E_0f\Vert^2 \right)
\end{equation*}
for all $f\in \bsob^1 (L(1),\mu_{\Sa})$ and all $0 < \eps \leq 1 - (\lambda_0/\lambda_1)$. (iii) For all $f\in \bsob^1 (L(1),\mu_{\Sa})$ and all $0 < \eps \leq 1 - (\lambda_0/\lambda_1)$, we have
\begin{equation*}
q_{\Sa,1}(f) \leq q_{\Sa,\eps}^0(f) + \lambda_0 \Vert E_0f\Vert^2.
\end{equation*}
\end{Lemma}

\noindent{\bf Proof.} (i) The first inequality follows from closedness of $q_{\Sa,1}$ and the second one for instance from \cite{Tay:99}, Prop. 5.2, p. 292. (ii) By the assumption on $\eps$, $\lambda_k \leq (\lambda_k - \lambda_0)/\eps$ for all $k\geq 1$. Hence,
\begin{eqnarray*}
q_{V}(f) &=&  \sum_{k\geq 0}\lambda_k \Vert E_kf\Vert^2 \leq \lambda_0\Vert E_0f\Vert^2 + \eps\sum_{k\geq 1}\frac{\lambda_k - \lambda_0}{\eps^2} \Vert E_kf\Vert^2 .
\end{eqnarray*}
By $q_H(f)\geq 0$ and $k_{\Sa}:=\max\lbrace 1,\lambda_0\rbrace$, we obtain the statement by
\begin{eqnarray*}
q_{V}(f) &\leq& \lambda_0\Vert E_0f\Vert^2 + \eps q_H(f) + \eps(q_{V}(f) - \lambda_0\Vert f\Vert^2)/\eps^ 2\\
&=& \eps q_{\Sa,\eps}^0(f)   +\lambda_0\Vert E_0f\Vert^2  \\
&\leq& k_{\Sa} \left( \eps q_{\Sa,\eps}^0(f) + \Vert E_0f\Vert^2 \right).
\end{eqnarray*}
(iii) By the spectral decomposition and the assumption on $\eps$ (in particular $\eps < 1$), we have, using again $\lambda_k\leq\eps^{-2}(\lambda_k-\lambda_0)$ for all $k\geq 1$,
\begin{eqnarray*}
q_{\Sa,1}(f) &=& q_V(f) + q_H(f) = \lambda_0 \Vert E_0f\Vert^2  + \sum_{k\geq 1}\lambda_k\Vert E_kf\Vert^2 + q_H(f)\\
&\leq& \lambda_0 \Vert E_0f\Vert^2  + \eps^{-1}\sum_{k\geq 1}\lambda_k\Vert E_kf\Vert^2  + q_H(f)\\
&\leq& \lambda_0 \Vert E_0f\Vert^2   + \sum_{k\geq 1}\eps^{-2}(\lambda_k-\lambda_0) \Vert E_kf\Vert^2 + q_H(f)\\
&=& q_{\Sa,\eps}^0 (f) + \lambda_0 \Vert E_0f\Vert^2.
\end{eqnarray*}
\hfill$\Box$

\subsection{Perturbation and Relative Boundedness} 

\noindent The following proposition summarizes all analytic facts that are needed for the analysis of the perturbation and that are consequences of the geometry of the tube and of the Dirichlet Laplacian. The proposition will be proved in Section \ref{geometry}.

\begin{Proposition}\label{pertprop}
1. The eigenspace $E_0$ of the vertical operator consists of functions of the form $f = \overline{f_b}\phi_0$, where $\overline{f_b}=f_b\circ\pi$ is a basic function and $\phi_0$ is constructed from the normalized eigenfunction $\varphi_0$ to the lowest eigenvalue $\lambda_0 > 0$ of the Dirichlet Laplacian on the flat unit ball by $\phi_0(x) = \varphi_0(d_{\Sa}(x,L))$. This is well defined since $\varphi_0$ is invariant with respect to orthogonal transformations. 
2. With the notations above, we have
$$
q_H(f) = \int_L \langle df_b,df_b\rangle_L \, d\mu_L
$$ 
for all $f\in E_0\cap\DD_H$.
3. For $0<\eps\leq 1$, the quadratic form
$$
l_{\eps}(f) := q_{\eps}(f)-q_{\Sa,\eps}(f) = \int_{L(1)} df \wedge (\#_{\eps} - \#_{\Sa,\eps}) df 
$$
with domain $\DD := \bsob^1(L(1),\mu_{\Sa})$ equals the sum
\begin{equation}
l_{\eps}(f) = \Omega (f) + r_{\eps}(f) 
\end{equation}
of two quadratic forms with the following properties:
\begin{itemize}
\item[{\rm (i)}] There is a constant $k_r > 0$, not depending on $\eps$, such that
$$
\left\vert r_{\eps}(f) \right\vert\leq \eps k_{r} \Vert f\Vert^2_{\bsob^1(L(1),\mu_{\Sa})}
$$
for all $0<\eps\leq 1$.
\item[{\rm (ii)}] There is some uniform constant $k_{\Omega} > 0$ such that 
$$
\left\vert \Omega(f)\right\vert \leq k_{\Omega} \,q_V(f) 
$$
for all $f\in \mathcal{D}_V$.
\item[{\rm (iii)}] The form $\Omega$ annihilates the eigenspace of $\Delta_V$ belonging to the smallest eigenvalue $\lambda_0$. Furthermore, if $E_0^{\perp}= 1-E_0$ denotes the projection onto the orthogonal complement of this eigenspace, we have $\Omega = \Omega\circ E_0^{\perp}$.
\end{itemize}
\end{Proposition}
As a consequence, the perturbation satisfies a Kato-type inequality with respect to $q_{\Sa,\eps}$. 

\begin{Corollary}\label{kato} There is a constant $k_l > 0$ such that 
\begin{equation*}
\left\vert l_{\eps}(f)\right\vert \leq  k_l\eps \left(q_{\Sa,\eps}^0(f) + \Vert f\Vert^2_{\bsob^1(L(1),\mu_{\Sa})} \right)
\end{equation*}
for all $f\in\bsob^1(L(1),\mu_{\Sa})$, $\eps < 1 - \lambda_0/\lambda_1$.
\end{Corollary}

\noindent{\bf Proof.} By Proposition \ref{pertprop}, 3(ii) and Lemma \ref{est1} (ii), we have
\begin{eqnarray*}
\left\vert \Omega (f)\right\vert &\leq& k_{\Omega} \,q_V(f)\\
&\leq& k_{\Omega}k_{\Sa}(\eps q^0_{\Sa,\eps}(f) + \Vert E_0f\Vert^2) .
\end{eqnarray*}
By Proposition \ref{pertprop}, 3(iii), we obtain
\begin{eqnarray*}
\left\vert \Omega (f)\right\vert &=& \left\vert \Omega (E_0^{\perp}f)\right\vert \leq k_{\Omega} k_{\Sa}\left(\eps\,  q_{\Sa,\eps}^0(E_0^{\perp}f) + \Vert E_0E_0^{\perp}f\Vert^2\right)\\
&=& \eps\,k_{\Omega}k_{\Sa}q_{\Sa,\eps}^0(E_0^{\perp}f).
\end{eqnarray*}
By \ref{sasaki_can_var}.c, 
\begin{eqnarray*}
q_H\left(\sum_{k\geq 0}E_kf\right)&=& \left\langle\sum_{k\geq 0}E_kf,\Delta_H \sum_{j\geq 0}E_jf\right\rangle = \sum_{k,j\geq 0}\langle f,\Delta_HE_kE_jf\rangle \\
&=& \sum_{k\geq 0}\langle E_k f,\Delta_H E_k f\rangle = \sum_{k\geq 0} q_H (E_k f)
\end{eqnarray*}
and $q_V(f)-\lambda_0\Vert f\Vert^2 = q_V(E_0^{\perp}f)-\lambda_0\Vert E_0^{\perp}f\Vert^2$, we obtain
$$
q_{\Sa,\eps}^0(E_0^{\perp}f) =\frac{q_V(E_0^{\perp}f)-\lambda_0\Vert E_0^{\perp}f\Vert^2}{\eps^2}+ q_H(E_0^{\perp}f)\leq q_{\Sa,\eps}^0(f)
$$ 
and hence
$$
\left\vert \Omega (f)\right\vert \leq    \eps\,k_{\Omega}k_{\Sa} q^0_{\Sa,\eps}(f) .
$$
By Proposition \ref{pertprop} (3i),
\begin{eqnarray*}
\left\vert l_{\eps}(f)\right\vert &\leq& \left\vert \Omega (f)\right\vert + \left\vert r_{\eps}(f)\right\vert \\
&\leq& \eps\,k_{\Omega} k_{\Sa} q_{\Sa,\eps}^0(f)  + \eps k_{r} \Vert f\Vert^2_{\bsob^1(L(1),\mu_{\Sa})} \\
&\leq& \eps \left(k_{\Omega}k_{\Sa} q_{\Sa,\eps}^0(f) +  k_r \Vert f\Vert^2_{\bsob^1(L(1),\mu_{\Sa})}\right).
\end{eqnarray*}
Letting $k_l := \max\lbrace k_{\Omega}k_{\Sa}, k_r \rbrace$ yields the statement.\hfill $\Box$ \\

\subsection{Equi-Coercivity and Convergence of the Minimizers}

Let $w\in L^2(L(1),\mu_{\Sa})$, $\alpha>0$ and $\phi_{\eps,\alpha,w}: \bsob^1 (L(1),\mu_{\Sa}) \to \R\cup\lbrace\infty\rbrace$ be the family given by
\begin{equation}\label{family}
\phi_{\eps,\alpha,w}(f) := \frac{1}{2}\left(q_{\eps}^{0}(f) + \alpha \Vert f \Vert^2\right) - \langle w,f\rangle 
\end{equation}
where $0 < \eps\leq 1$. Recall that $\Vert - \Vert$ and $\langle -,-\rangle$ denote the norm and the scalar product on $L^2(L(1),\mu_{\Sa})$. The Kato-type inequality Corollary \ref{kato} implies the following fundamental result for the perturbation family:

\begin{Proposition}\label{coerz} Let $\lambda_0 > 0$ denote the eigenvalue from \ref{sasaki_can_var}.b. For all $\alpha > \lambda_0$, we have with the constant $A > 0$ from Lemma \ref{est1} (i) 
\begin{equation*}
q_{\eps}^0 (f) + \alpha \Vert f\Vert^2 \geq \frac{1}{2A}\Vert f\Vert_{\bsob^1 (L(1),\mu_{\Sa})}^2
\end{equation*} 
for all $0<\eps < \eps_0 := \min\lbrace (2 k_l (1+A))^{-1}, 1- \lambda_0\lambda_1^{-1}\rbrace$.
\end{Proposition}

\noindent{\bf Proof.} By Corollary \ref{kato} and the assumptions on $\eps$
\begin{eqnarray*}
q^0_{\eps}(f) &=& q^0_{\Sa,\eps}(f) + l_{\eps}(f) \geq q^0_{\Sa,\eps}(f) - \left\vert l_{\eps}(f) \right\vert \\
&\geq& q^0_{\Sa,\eps}(f) - \eps k_l q^0_{\Sa,\eps}(f) -  \eps k_l \Vert f \Vert^2_{\bsob^1(L(1),\mu_{\Sa})} \\
&=& (1-\eps k_l) q^0_{\Sa,\eps}(f) - \eps k_l \Vert f\Vert^2_{\bsob^1(L(1),\mu_{\Sa})}.
\end{eqnarray*}
By Lemma \ref{est1} (i) and (iii), we have
$$
\frac{1}{A}\Vert f\Vert^2_{\bsob^1(L(1),\mu_{\Sa})}\leq q_{\Sa,1}(f)\leq q^0_{\Sa,\eps}(f) + \lambda_0\Vert E_0f\Vert^2.
$$
That implies
\begin{eqnarray*}
q^0_{\eps}(f) &\geq& (1-\eps k_l) ( \frac{1}{A}\Vert f\Vert^2_{\bsob^1(L(1),\mu_{\Sa})}- \lambda_0\Vert E_0f\Vert^2 ) - \eps k_l \Vert f\Vert^2_{\bsob^1(L(1),\mu_{\Sa})} \\
&=& \frac{1-\eps k_l(A+1)}{A} \Vert f \Vert^2_{\bsob^1(L(1),\mu_{\Sa})} - (1-\eps k_l)  \lambda_0\Vert E_0f\Vert^2  ,
\end{eqnarray*}
hence for $\alpha > \lambda_0$
\begin{eqnarray*}
	q^0_{\eps}(f) + \alpha\Vert f\Vert^2 &\geq& \frac{1-\eps k_l(A+1)}{A} \Vert f \Vert^2_{\bsob^1(L(1),\mu_{\Sa})} \geq \frac{1}{2A}\Vert f \Vert^2_{\bsob^1(L(1),\mu_{\Sa})}
\end{eqnarray*}
for $\eps\leq (2 k_l (1+A))^{-1}$.\hfill$\Box$\\

\noindent From now on, we will always denote the parameter bound by
\begin{equation*}
\eps_0 := \min\lbrace (2 k_l (1+A))^{-1}, 1- \lambda_0\lambda_1^{-1}\rbrace .
\end{equation*}
We now draw some conclusions concerning the family (\ref{family}).

\noindent a. As a first consequence, the functions are lower semi-continuous with respect to the weak topology on the boundary Sobolev space $\bsob^1(L(1),\mu_{\Sa})$.

\begin{Corollary} Let $\alpha > \lambda_0$, $w\in L^2(L(1),\mu_{\Sa})$. Then, the functions $\phi_{\eps,\alpha,w}:\bsob^1(L(1),\mu_{\Sa})\to \R$ are continuous in the strong, and lower semi - continuous in the weak topology on the boundary Sobolev space for all $0<\eps < \eps_0$.  
\end{Corollary}

\noindent{\bf Proof.} (i) Continuity in the strong topology follows from closedness of $q_{\eps}^0$. (ii) By Proposition \ref{coerz}, $\phi_{\eps,\alpha,0}(f)\geq 0$ for all $f\in\bsob^1(L(1),\mu_{\Sa})$, i.e. the quadratic form $\phi_{\eps,\alpha,0}$ is non-negative. That implies for $f_n,f\in\bsob^1(L(1),\mu_{\Sa})$
\begin{eqnarray*}
0 &\leq& \phi_{\eps,\alpha,0} (f_n-f) \\
&=&  \phi_{\eps,\alpha,0} (f_n) + \phi_{\eps,\alpha,0} (f) -2 \left( b_{\eps}(f,f_n) + (\alpha-\lambda_0 \eps^{-2})\langle f_n,f\rangle\right) .
\end{eqnarray*}
For $f_n\rightharpoonup f$ weakly, we have $\lim_{n\to\infty}b_{\eps}(f,f_n) =b_{\eps}(f,f)$. That implies $\lim_{n\to\infty}\langle f_n,f\rangle = \Vert f\Vert^2$ and therefore, $0 \leq \liminf_n  \phi_{\eps,\alpha,0} (f_n) - \phi_{\eps,\alpha,0} (f)$.
By $\lim_{n\to\infty} \langle w,f_n\rangle = \langle w,f\rangle$, we obtain $\liminf_n  \phi_{\eps,\alpha,w} (f_n) = \phi_{\eps,\alpha,w} (f)$.
\hfill$\Box$ \\

\noindent b. The second consequence of the estimate is the following uniform statement about the location of the spectrum.

\begin{Corollary}\label{spc} The operator
\begin{equation*}
H_{\eps}^0 + \lambda_0 =  H_{\eps}  -\frac{\lambda_0}{\eps^2} + \lambda_0\geq 0
\end{equation*}
is non-negative uniformly for all $0<\eps < \eps_0$. In particular, the operator is self-adjoint with $\mathrm{spec} (H_{\eps}^0 + \alpha)\subset \lbrack \alpha-\lambda_0,\infty)$ for all $\alpha > \lambda_0$.
\end{Corollary}

\noindent The functions $\phi_{\eps,\alpha,w}$ are strictly convex and differentiable with differential $\nabla \phi_{\eps,\alpha,w} : \bsob^1 (L(1),\mu_{\Sa}) \to\bsob^1 (L(1),\mu_{\Sa})^*$. The minimizer $f^*_{\eps,\alpha,w}$ is therefore unique and satisfies  
\begin{equation*}
\nabla_{f^*_{\eps,\alpha,w}} \phi_{\eps,\alpha,w} (u) = 0 
\end{equation*}
for all $u\in\bsob^1 (L(1),\mu_{\Sa})$. This is equivalent to the statement that $f^*_{\eps,\alpha,w}$ is a weak solution of 
\begin{equation*}
(H_{\eps}^0 + \alpha)f^*_{\eps,\alpha,w} = w.
\end{equation*}
However, by Corollary \ref{spc}, there is indeed a strong solution given by the resolvent 
\begin{equation*}
f^*_{\eps,\alpha,w} = (H_{\eps}^0 + \alpha)^{-1} w\in \bsob^1\cap\sob^2 (L(1),\mu_{\Sa}).
\end{equation*}

\begin{Corollary} Let $\alpha > \lambda_0 +1$ and $w\in L^2(L(1),\mu_{\Sa})$ be fixed. Then, the set
\begin{equation*}
\EE (t) := \bigcup_{0<\eps < \eps_0} \left\lbrace f\in\bsob^1 (L(1),\mu_{\Sa})\, :\, \phi_{\eps,\alpha,w}(f) \leq t\right\rbrace \subset \bsob^1(L(1),\mu_{\Sa})
\end{equation*}
is norm bounded for all $t\in\R$.
\end{Corollary}

\noindent{\bf Proof.} By Proposition \ref{coerz}
\begin{eqnarray*}
\frac{1}{4A}\Vert f\Vert_{\bsob^1 (L(1),\mu_{\Sa})}^2 &\leq& \phi_{\eps,\alpha-1,0}(f) = \phi_{\eps,\alpha-1,w}(f) + \langle w,f\rangle \\
&\leq& \phi_{\eps,\alpha-1,w}(f) + \frac{1}{2}\left(\Vert f\Vert^2 + \Vert w\Vert^2\right) \\
&=& \phi_{\eps,\alpha,w}(f) + \frac{1}{2}\Vert w\Vert^2 \\
\end{eqnarray*}
uniformly for all $0<\eps<\eps_0$. Hence, $f\in \EE (t)$ implies $\Vert f\Vert_{\bsob^1 (L(1),\mu_{\Sa})}^2\leq 2A(2t + \Vert w \Vert^2)$ and $\EE (t)$ is indeed norm - bounded.\hfill$\Box$\\

\noindent Equi-coercivity implies that every sequence of minimizers contains a convergent subsequence in the following sense.

\begin{Corollary}\label{convergence} Let $\eps_n > 0$ such that $\lim_{n\to\infty} \eps_n = 0$ and $\alpha > \lambda_0 + 1$. Denote by 
$$
f_n := f^*_{\eps_n,\alpha,w} = \left( H_{\eps_n}^0 + \alpha\right)^{-1} w
$$
the sequence of the (unique) minimizers of the functionals $\phi_{\eps_n,\alpha,w}$. Then, $f_n$ contains a subsequence which converges strongly in $L^2 (L(1),\mu_{\Sa})$ and weakly in $\bsob^1 (L(1),\mu_{\Sa})$.
\end{Corollary}

\noindent{\bf Proof.} Let $u\in E_0$. Then for $K := \sup_{0 < \eps < \eps_0} \eps r_{ \eps}(u)$, we have by Proposition \ref{pertprop}, (1) and (3 iii)
$$
\phi_{\eps_n,\alpha,w}(u) = q^0_{\Sa,\eps_n}(u) + l_{\eps_n} (u) \leq \frac{1}{2} q_H(u) + \alpha \Vert u\Vert^2 - \langle u,w\rangle + K =: t_0
$$
independent of $n\geq 1$. Hence, for all $n\geq 1$, we have $\phi_{\eps_n,\alpha,w}( f^*_{\eps_n,\alpha,w})\leq t_0$ and the sequence is therefore contained in a subset which is norm-bounded in $\bsob^1 (L(1),\mu_{\Sa})$. Hence, the subset is weakly compact in $\bsob^1 (L(1),\mu_{\Sa})$ and compact in $L^2 (L(1),\mu_{\Sa})$.\hfill$\Box$\\

\subsection{Epi-Convergence and Convergence of Semigroups in $\mathbf{L^2(L(1),\mu_{\Sa})}$}

By Corollary \ref{convergence} above, every sequence $(f_n)$ contains a convergent subsequence. Now we are going to identify the limit, which will also prove that the sequence $( H_{\eps_n}^0 + \alpha)^{-1}$, $n\geq 1$, of resolvent operators converges strongly in the space of bounded operators on the Hilbert space $L^2(L(1),\mu_{\Sa})$. Recall from Proposition \ref{pertprop} that we may decompose a given $f\in E_0$ by $f = \overline{f_b}\phi_0$, where the function $f_b$ on $L$ is almost surely given by $f_b(x) := \langle f_x,\phi_0\rangle_{\pi^{-1}(x)}$, where $\langle -,-\rangle_{\pi^{-1}(x)}$ denotes the scalar product on the fibre with the metric induced by the Sasaki metric, and $\phi_0$ only depends on the radial distance to the submanifold.

\begin{Proposition}\label{epi} For $\alpha > \lambda_0 + 1$, the functions $\phi_{\eps,\alpha,w}$ epi-converge, as $\eps$ tends to zero, to
\begin{equation*}
\phi_{\alpha,w}(f) := \left\lbrace \begin{array}{ll}
                                    \frac{1}{2}\int_L  \left(\langle df_b, df_b \rangle_L + \alpha f_b^2\right)\,d\mu_L - \langle f,w \rangle & ,f\in E_0 \\
                                    \infty & ,f\notin E_0
                                   \end{array}
\right. 
\end{equation*}
with respect to the weak topology on $\bsob^1(L(1),\mu_{\Sa})$.
\end{Proposition}

\noindent{\bf Proof.} For $\alpha \geq \lambda_0$, the function $\phi_{\eps,\alpha,0}$ is a non - negative bounded quadratic form and therefore weakly lower semi - continuous on the boundary Sobolev space $\bsob^1(L(1),\mu_{\Sa})$. We prove the assertion in three steps: 

\noindent a. Let $(\eps_n)_{n\geq 1}$ be a decreasing sequence which converges to zero. The functions 
$$
q_{n}(f) = (2\eps_n)^{-1} \left( q_V (f) - \lambda_0 \Vert f\Vert^2\right) 
$$
form an increasing family $q_n : \bsob^1(L(1),\mu_{\Sa})\to \R$ of non - negative quadratic forms. Thus, the $q_n$ are lower semi - continuous in the weak topology, and, by \cite{Dal:93}, Proposition 5.4, p. 47, epi - converge to
$$
q_{\infty}(f) := \left\lbrace \begin{array}{ll}
                                    \sup_n q_n(f) & \mathrm{,if \, the \, limit \, exists} \\
                                    \infty & \mathrm{,else}
                                   \end{array}
\right. = \left\lbrace \begin{array}{ll}
                                    0 & ,f\in E_0 \\
                                    \infty & ,f\notin E_0
                                   \end{array}
\right. .
$$
b. For the Sasaki metric, we have 
$$
q'(f) := \phi_{\Sa, \eps_n, \alpha,w}(f) - q_n(f) = \frac{1}{2}\left(q_H(f) + \alpha\Vert f\Vert^2\right)-\langle w,f\rangle.
$$
By $\alpha \geq 0$, the right hand side $2^{-1}(q_H(f) + \alpha\Vert f\Vert^2)$ is a bounded and non-negative quadratic form on $\bsob^1(L(1),\mu_{\Sa})$ and therefore weakly lower semi-continuous. By the monotonicity in part a. that implies that $\phi_{\Sa, \eps_n, \alpha,w}(f)$ epi - converges to  
$$
\phi_{\Sa, \alpha,w}(f) = q_{\infty}(f) + q'(f) = \left\lbrace \begin{array}{ll}
                                    \frac{1}{2}q_H(f) + \alpha\Vert f\Vert^2-\langle w,f\rangle & ,f\in E_0 \\
                                    \infty & ,f\notin E_0
                                   \end{array}
\right. 
$$
in the weak topology (cf. \cite{Dal:93}, Example 6.24 (b), p. 64). c. For a general metric, we have by Proposition \ref{pertprop}
$$
q''(f) := \phi_{\eps_n, \alpha,w}(f) - \phi_{\Sa, \eps_n, \alpha,w}(f) = l_{\eps_n} (f),
$$
with the Kato type estimate for $l_{\eps}$ from Corollary \ref{kato}. Let $f\in\bsob^1(L(1),\mu_{\Sa})$ be fixed. By b., we have 
$$
\lim_{\eps\to 0} \phi_{\Sa, \eps_n, \alpha,w}(f) = \left\lbrace \begin{array}{ll}
                                    \frac{1}{2}(q_H(f) + \alpha\Vert f\Vert^2)-\langle w,f\rangle & ,f\in E_0 \\
                                    \infty & ,f\notin E_0
                                   \end{array}\right. .
$$
By Proposition \ref{pertprop}, (3), we have
$$
\lim_{\eps\to 0}l_{\eps} (f) = \Omega (f) + \lim_{\eps \to 0}  \eps r_{\eps}(f) = \Omega (f)
$$
and $\Omega=\Omega\circ E_0^{\perp}$. Thus 
\begin{equation}\label{epilimsup}
\lim_{\eps\to 0}\phi_{\eps_n, \alpha,w}(f) = \left\lbrace \begin{array}{ll}
                                    \frac{1}{2}(q_H(f) + \alpha\Vert f\Vert^2)-\langle w,f\rangle & ,f\in E_0 \\
                                    \infty & ,f\notin E_0
                                   \end{array}
\right..
\end{equation}
\noindent On the other hand,  we have by Corollary \ref{kato} for all $f\in\bsob^1 (L(1),\mu_{\Sa})$
\begin{eqnarray*}
& & \phi_{\eps_n, \alpha,w}(f) \\ 
&=& \phi_{\Sa, \eps_n, \alpha,w}(f) + \frac{1}{2}l_{\eps_n}(f)\geq  \phi_{\Sa, \eps_n, \alpha,w}(f) - \frac{1}{2}\vert  l_{\eps_n}(f )\vert \\
&=& \frac{1}{2}\left(q_{\Sa,\eps}^{0}(f) + \alpha \Vert f \Vert^2\right) - \langle w,f\rangle -\frac{\eps k_l}{2} (q_{\Sa,\eps}^{0}(f) + \Vert f\Vert_{\bsob^1 (L(1),\mu_{\Sa})}^2)\\
&=& \left\lbrack 1- \frac{\eps k_l}{2}\right\rbrack\phi_{\Sa, \eps_n, \alpha,w}(f) - \frac{\eps k_l}{2} \left\lbrack \alpha \Vert f \Vert^2  - \langle w,f\rangle  + \Vert f\Vert^2_{\bsob^1 (L(1),\mu_{\Sa})}\right\rbrack .
\end{eqnarray*}
Let now $(f_n)$ be a weakly convergent sequence in $\bsob^1 (L(1),\mu_{\Sa})$ with weak limit $f$. Since weakly convergent sequences are norm-bounded, we have 
$$
\sup \lbrace \Vert f_n\Vert^2, \vert\langle w,f_n\rangle\vert, \Vert f_n\Vert_{\bsob^1 (L(1),\mu_{\Sa})}^2\,:\, n\geq 1\rbrace =:M<\infty.
$$
Hence,
$$
\limsup_{n\to\infty} \frac{\eps_n k_l}{2} ( \alpha \Vert f_n \Vert^2  - \langle w,f_n\rangle  + \Vert f_n\Vert_{\bsob^1 (L(1),\mu_{\Sa})}^2)\leq \lim_{n\to\infty} \frac{3M\eps_n k_l}{2} = 0.
$$
and therefore,
\begin{equation}\label{epiliminf}
\liminf_{n\to\infty}  \phi_{\eps_n, \alpha,w}(f) \geq \liminf_{n\to\infty}  \phi_{\Sa,\eps_n, \alpha,w}(f) .
\end{equation}
(\ref{epilimsup}) and (\ref{epiliminf}) imply that the functionals associated to the induced metric epi-converge to the same limit functional as the functionals associated to the Sasaki-metric (cf. \cite{Dal:93}, Prop. 8.1, p. 87). By rewriting the limit using Proposition \ref{pertprop} (1),(2), we obtain the statement.\hfill$\Box$\\ 

\noindent Under epi-convergence of functionals, every convergent sequence of minimizers converges to a minimizer of the limit. That implies: 

\begin{Corollary} Let $\eps_n > 0$ such that $\lim_{n\to\infty} \eps_n = 0$ and $\alpha > \lambda_0 +1$. Denote by 
$$
f_n :=  \left( H_{\eps_n}^0 + \alpha\right)^{-1} w
$$
the sequence of the (unique) minimizers of the functionals $\phi_{\eps_n,\alpha,w}$. Then, 
$$
\lim_{n\to\infty} f_n = E_0\,\left( \Delta_L + \alpha\right)^{-1}\,E_0 w
$$
strongly in $L^2 (L(1),\mu_{\Sa})$ and weakly in $\bsob^1 (L(1),\mu_{\Sa})$.
\end{Corollary}

\noindent{\bf Proof.} By Corollary \ref{convergence}, every subsequence of $f_n$ contains a convergent subsequence. By Proposition \ref{epi} above, the epi-limit of the functionals $\phi_{\eps_n,\alpha,w}$ is given by $\phi_{\alpha,w}$. By \cite{Dal:93}, Corollary 7.20, p. 81, every convergent sequence of minimzers $f_n$ will converge to a minimizer of the limit functional weakly in $\bsob^1(L(1),\mu_{\Sa})$ and therefore also strongly in $L^2(L(1),\mu_{\Sa})$. However, the minimizer of $\phi_{\alpha,w}$ is unique and given by
$f_{\infty} = E_0\,\left( \Delta_L + \alpha\right)^{-1}\,E_0 w$. That implies the statement.\hfill$\Box$\\

\noindent Finally, strong convergence of the resolvents and uniform sectoriality of the corresponding operators imply convergence of the associated semigroups, even if the limit is just a pseudo-resolvent.

\begin{Proposition}\label{semi1} Let $f\in L^2(L(1),\mu_{\Sa})$. Then
$$
\lim_{\eps\to 0} e^{-\frac{t}{2}H_{\eps}^0}f = E_0\,e^{-\frac{t}{2}\Delta_L}\,E_0f.
$$
in $L^2(L(1),\mu_{\Sa})$.
\end{Proposition}

\noindent{\bf Proof.} Let $\alpha > \lambda_0 +1$. Then, by Corollary \ref{spc}, $\spec (H_{\eps}^0 + \alpha)\subset \lbrack 1,\infty)$ uniformly for all $\eps_0 > \eps > 0$. Thus, integration along a suitable contour $\gamma$ (for instance, the negatively oriented boundary of a sector $\lbrace z\in\C\,:\,\arg(z) < \pi/4\rbrace$) yields for $t>0$
$$
e^{-\frac{t}{2}(H_{\eps}^0+\alpha)}f = \frac{1}{2\pi i}\int_{\gamma} e^{-\frac{z t}{2}} R(H_{\eps}^0+\alpha,z)f \,dz
$$
and
$$
E_0\,e^{-\frac{t}{2}(\Delta_L + \alpha)}\,E_0f = \int_{\gamma} e^{-\frac{zt}{2}} E_0\,R(\Delta_L+\alpha,z)\,E_0f \, dz .
$$
Hence, the convergence result for the resolvents implies by dominated convergence  
$$
\lim_{\eps\to 0} e^{-\frac{t}{2}(H_{\eps}^0+\alpha)}f = E_0\,e^{-\frac{t}{2}(\Delta_L + \alpha)}\,E_0f.
$$
Multiplication by $e^{\alpha t/2}$ yields the statement.
\hfill$\Box$\\

\section{The Tube Geometry and a Proof of Proposition \ref{pertprop}}\label{geometry}

In this subsection, we investigate the local geometry of the tubes and its consequence for the behaviour of the asymptotic Dirichlet problem. \noindent In particular, we derive an asymptotic formula for the quadratic form (\ref{scalequad}). As a first result, we express the metric $g$ in terms of the Sasaki metric.

\subsection{Jacobi fields and the metric on the cotangent bundle} First of all, we have to fix some notation.

\begin{Definition} Let $W\in NL$ with $\pi(W)=x\in L$. {\rm (i)} By $\mathcal{L}_W : T_xL\oplus N_xL \to T_WNL$ we denote the {\em lift}
$$
\mathcal{L}_W (X,V) := X^h(W) + J_WV,
$$
where $X^h$ denotes the horizontal and $J_WV$ the vertical lift of the respective vectors. {\rm (ii)} By $\exp: NL \to M$, we denote $\exp (W) := \exp_{\pi (W)}(W)$, i.e. $\gamma_W(1)$, where $\gamma_W$ is the geodesic in $M$ with $\gamma_W(0)=x$, $\dot{\gamma}_{W}(0)= W$. {\rm (iii)} By $P_W:T_xNL\to T_WNL$ we denote parallel translation along $\gamma_W$.  
\end{Definition}

\noindent Note that $\mathcal{L}_W$ depends only on the geometry of the normal bundle, i.e. the induced connection on $NL$, whereas $\exp$ and $P$ depend on the geometry of the ambient space. We formulate the essence of what we need from the theory of Jacobi fields in the following way:

\begin{Proposition}\label{jacobi_fields} {\rm (i)} There is an endomorphism $\mathcal{A}_W:T_xNL \to T_xNL$ such that
\begin{equation}\label{umlauf}
\mathcal{A}_W = P_W^{-1}\circ\exp_{*,W}\circ\,\mathcal{L}_W.
\end{equation}
{\rm (ii)} For $X\in T_xL\subset T_xNL$, we have
\begin{equation}\label{horizontal_umlauf}
\mathcal{A}_W X = X - A_WX + \frac{1}{2}R(W,X)W + O(\Vert W\Vert^3),
\end{equation}
where $A_W$ is the Weingarten map of the embedding $L\subset M$ and $R$ is the curvature tensor of $M$ at $x$. {\rm (iii)} For $V\in N_xL\subset T_xNL$, we have
\begin{equation}\label{vertikal_umlauf}
\mathcal{A}_W V = V + \frac{1}{6}R(W,V)W + O(\Vert W\Vert^3).
\end{equation}
\end{Proposition}

\noindent{\bf Proof.} This is a standard Jacobi field argument, cf. \cite{Cha:93}, p. 132 ff., for a somewhat different formulation.
\hfill$\Box$ \\

\noindent{\bf Remark.} In particular, if $M=NL$ is equipped with the Sasaki metric, we have $\mathcal{L}_W = P_W$, $\exp_{*,W}= \mathrm{id}_{T_WNL}$ and therefore, $\mathcal{A}_{W}^{\Sa} = \mathrm{id}_{T_xNL}$. \\

\begin{Lemma}\label{metric1} Let $W\in NL$ and $\langle-,-\rangle$ denote the pullback of the metric on $L(1)$ via $\exp$. Let $\langle -,-\rangle_{\Sa}$ denote the Sasaki metric on $NL$ and $\mathcal{A}_W$ the map from Proposition \ref{jacobi_fields}, (i). Then, we have for all $Z,Z'\in T_WNL$ 
\begin{equation*}
\langle  Z,  Z'\rangle_W = \langle U_W Z, U_W Z'\rangle_{\Sa,W}
\end{equation*} 
where $U_W: T_WNL \to T_WNL$ is given by $U_W := \mathcal{L}_W \mathcal{A}_W \mathcal{L}_W^{-1}$.
\end{Lemma}

\noindent{\bf Proof.} By definition, $\exp$ is an isometry. Hence, by the remark above
\begin{eqnarray*}
\langle  Z,  Z'\rangle_W &=& \langle \mathcal{L}_W\mathcal{L}_W^{-1} Z, \mathcal{L}_W \mathcal{L}_W^{-1}Z'\rangle_W \\
&=& \langle \exp_{*,W}\mathcal{L}_W \mathcal{L}_W^{-1}Z, \exp_{*,W}\mathcal{L}_W \mathcal{L}_W^{-1}Z'\rangle_{\gamma_W(1)} \\
&=& \langle P_W^{-1}\exp_{*,W}\mathcal{L}_W \mathcal{L}_W^{-1}Z, P_W^{-1}\exp_{*,W}\mathcal{L}_W \mathcal{L}_W^{-1}Z'\rangle_{T_xL \oplus N_xL} \\
&=& \langle \mathcal{A}_W\mathcal{L}_W^{-1} Z, \mathcal{A}_W \mathcal{L}_W^{-1}Z'\rangle_{T_xL\oplus N_xL} \\
&=&\langle \mathcal{L}_W \mathcal{A}_W \mathcal{L}_W^{-1}Z, \mathcal{L}_W \mathcal{A}_W\mathcal{L}_W^{-1} Z'\rangle_{\Sa,W}.
\end{eqnarray*}
\hfill$\Box$ \\

\noindent{\bf Remark.} By $\LL^{-1}=\pi_*\oplus K$, and the definition of the Sasaki metric, we have $\langle  Z,  Z'\rangle_W = \langle \mathcal{A}_W \mathcal{L}_W^{-1}Z,  \mathcal{A}_W \mathcal{L}_W^{-1}Z'\rangle_{L,x} + \langle \mathcal{A}_W \mathcal{L}_W^{-1}Z,  \mathcal{A}_W \mathcal{L}_W^{-1}Z'\rangle_{N_{x}L}.$\\

\noindent To describe the effect of the rescaling on the dual metric, we first decompose the cotangent bundle similarly to the tangent bundle. Therefore, we note that the dual maps $\LL_W^*: T^*_WNL\to T^*_xL\oplus N^*_xL$ and $\left(\LL_W^{-1}\right)^*:T^*_xL\oplus N^*_xL \to T^*_WNL$ are given by $\LL_W^*\eta = (\kappa_W\eta, J_W^*\eta)$ and $\LL_W^{-1 \,*}(\xi,\omega) = \pi^*\xi + K_W^*\omega$, where $\kappa_W$ denotes the dual of the horizontal lift. That implies that $\eta\in T_W^*NL$ can be uniquely written as $\eta = \LL_W^{-1\,*}(\eta^{\top},\eta^{\perp})$, with $\eta^{\top}\in T^*_xL$, $\eta^{\perp}\in N^*_xL$. \\

\noindent Consider the orthogonal decomposition $T_xNL = T_xL \oplus N_xL$ with orthogonal projections $P_T:T_xNL \to T_xL, P_N:T_xNL \to N_xL$. Thus, Proposition \ref{jacobi_fields} reads
$$
\mathcal{A}_{\eps W}  = \mathrm{id} - \eps A_W P_T + \eps^2R_W (2^{-1}P_T + 6^{-1} P_N)  + O( \eps^3),
$$
with $R_WZ = R(W,Z)W$. That implies 

\begin{Lemma}\label{dual_A} $\mathcal{A}_{\eps W}^{-1\,*}:T^*_{x}NL \to T^*_{x}NL$ is given by
$$
\mathcal{A}_{\eps W}^{-1\,*} = \mathrm{id} + \eps A_W^* P_T^* - \eps^2\left((2^{-1} P_T^* + 6^{-1}P_N^*)R_W^* -  P_T^*A_W^{*\,2} \right) + O( \eps^3).
$$
\end{Lemma}

\noindent{\bf Proof.} By
$$
\mathcal{A}_{\eps W}^{-1} =\mathrm{id} + \eps  A_W P_T - \eps^2\left( R_W (2^{-1} P_T + 6^{-1} P_N)-  A_W^2P_T \right) + O( \eps^3),
$$
we have with $Z=(X,V)\in T_xNL$
\begin{eqnarray*}
\mathcal{A}_{\eps W}^{-1\,*}\eta^{\top}(Z) &=& \eta^{\top}(X + \eps A_WX - \eps^2(2^{-1}R_W-A_W^{2})X - 6^{-1}R_WV)\\
&&+O(\eps^3) \\
&=& \eta^{\top}(X) + \eps \eta^{\top}(A_WX) - \eps^2\eta^{\top}(R_W(2^{-1}X + 6^{-1}V))\\
& & +\eps^2\eta^{\top}(A_W^{2}X) +O(\eps^3).
\end{eqnarray*}
The {\em Weingarten map} $A_W$ is an endomorphism of the tangent space $TL$. Hence, $A_WP_T=P_TA_WP_T$, $A_W^2P_T=P_TA_W^2P_T$ and we have
\begin{eqnarray*}
\mathcal{A}_{\eps W}^{-1\,*}\eta^{\perp}(Z) &=& \eta^{\perp}(V  - \eps^2R_W(2^{-1}X + 6^{-1}V)+O(\eps^3)) \\
&=& \eta^\perp(V)  - \eps^2\eta^{\perp}(R_W(2^{-1}X + 6^{-1}V))+O(\eps^3) .
\end{eqnarray*}
That implies the statement.\hfill$\Box$ \\

\begin{Lemma}\label{dual_B} Let $\sigma_{\eps}^*:T^*_{W}NL \to T^*_{\eps W}NL$ and $\eta = \pi^*\eta^{\top} + K_W^*\eta^{\perp}\in T_W^*NL$. Then,
$$
\LL_{\eps W}^{*}\sigma_{\eps}^*\eta = (P_T^* + \eps^{-1}P_N^*)\LL_{W}^{*}\eta = (\eta^{\top}, \eps^{-1}\eta^{\perp}).
$$
\end{Lemma}

\noindent{\bf Proof.} By $\pi\circ\sigma_{\eps} = \pi$, we have 
$$
\LL_{\eps W}^*\sigma_{\eps}^*\pi^*(\eta^{\top}) = \LL_{\eps W}^*(\pi^*\eta^{\top}) = \eta^{\top},
$$ 
and by $K_W\circ J_W = \mathrm{id}_{NL}$ together with $\sigma_{\eps\,*}J_{\eps W} = \eps^{-1}J_W$, we obtain 
$$
\LL_{\eps W}^*\sigma_{\eps}^*K_W^*\eta^{\perp}= J_{\eps W}^*\sigma_{\eps}^{*}K_W^*\eta^{\perp} =\eps^{-1}J_W^*K_{ W}^*\eta^{\perp}=\eps^{-1}\eta^{\perp}.$$
\hfill$\Box$ \\

\noindent Now we come to the first statement, a representation of the induced metric in a small tubular neighbourhood around the submanifold. 

\begin{Proposition}\label{induced_dual} Let $W\in L(1)$. The induced rescaled metric $\langle -,-\rangle_{\eps,W}$ on the cotangent bundle $T^*L(1)$ is asymptotically given by 
$$
\langle \eta,\zeta\rangle_{\eps,W} = \langle \eta,\zeta\rangle_{\Sa,\eps,W} - \frac{1}{3}\langle \eta^{\perp},R_W^*\zeta^{\perp}\rangle_{\Sa,\pi (W)} + \eps r_{\eps}(\eta,\zeta).
$$
Here, $r_{\eps}:T^*L(1)\times T^*L(1)\to\R$ denotes a bilinear form with smooth coefficients that are uniformly bounded in $0<\eps\leq 1$.
\end{Proposition}

\noindent{\bf Proof.} Let $\eta = \LL_W^{-1\,*}(\eta^{\top},\eta^{\perp})$ and $\zeta = \LL_W^{-1\,*}(\zeta^{\top},\zeta^{\perp})$. By Lemma \ref{dual_A} and Lemma \ref{dual_B}
\begin{eqnarray*}
\langle \eta,\zeta \rangle_{\eps, W} &:=& \langle \sigma_{\eps}^*\eta,\sigma_{\eps}^*\zeta \rangle_{\eps W} = \langle U_{\eps W}^{-1\,*}\sigma_{\eps}^*\eta,U_{\eps W}^{-1\,*}\sigma_{\eps}^*\zeta \rangle_{\Sa,\eps W}\\
&=& \langle \LL_{\eps W}^*U_{\eps W}^{-1\,*}\sigma_{\eps}^*\eta, \LL_{\eps W}^*U_{\eps W}^{-1\,*}\sigma_{\eps}^*\zeta \rangle_{\Sa,x} \\
&=& \langle \mathcal{A}_{\eps W}^{-1\,*}\LL_{\eps W}^{*}\sigma_{\eps}^*\eta, \mathcal{A}_{\eps W}^{-1\,*}\LL_{\eps W}^*\sigma_{\eps}^*\zeta \rangle_{\Sa,x}\\
&=& \langle \mathcal{A}_{\eps W}^{-1\,*}(\eta^{\top},\eps^{-1}\eta^{\perp}), \mathcal{A}_{\eps W}^{-1\,*}(\zeta^{\top},\eps^{-1}\zeta^{\perp}) \rangle_{\Sa,x}.
\end{eqnarray*}
By Lemma \ref{dual_A}, we obtain
\begin{eqnarray*}
\mathcal{A}_{\eps W}^{-1\,*}(\eta^{\top},0) &=& ( \eta^{\top} + \eps (A_W^*\eta^{\top} - 6^{-1}(R_W^*\eta^\perp)^{\top}), 0) + O(\eps^2) \\
\mathcal{A}_{\eps W}^{-1\,*}(0,\eps^{-1}\eta^{\perp}) &=& ( -\eps 6^{-1}(R_W^*\eta^\perp)^{\top},\frac{1}{\eps}\eta^{\perp} - \eps 6^{-1}(R_W^*\eta^\perp)^{\perp}) + O(\eps^2) .
\end{eqnarray*}
That implies by $P_NP_T=P_TP_N=0$ and by the symmetries of the curvature tensor
$$
\langle \eta,\zeta \rangle_{\eps, W} 
= \langle \eta^{\top},\zeta^{\top}\rangle_{\Sa,x} + \frac{1}{\eps^2}\langle \eta^{\perp},\zeta^{\perp}\rangle_{\Sa,x} -\frac{1}{3}\langle \eta^{\perp},R_W^*\zeta^{\perp}\rangle_{\Sa,x} + O(\eps) .
$$
Finally, by Lemma \ref{dual_B} and by the definition of the Sasaki metric 
\begin{eqnarray*}
\langle \eta^{\top},\zeta^{\top}\rangle_{\Sa,x} + \frac{1}{\eps^2}\langle \eta^{\perp},\zeta^{\perp}\rangle_{\Sa,x} &=& \langle \LL_{\eps W}^{*}\sigma_{\eps}^*\eta,\LL_{\eps W}^{*}\sigma_{\eps}^*\zeta\rangle_{\Sa,x} \\
&=& \langle \sigma_{\eps}^*\eta,\sigma_{\eps}^*\zeta\rangle_{\Sa,\eps W}.
\end{eqnarray*}
\hfill $\Box$\\

\noindent That means, the leading term in the expansion of the canonical variation of the induced metric is given by the canonical variation of the Sasaki metric. Furthermore, there is only one additional term of relevant order given by a curvature form on the fibres. 

\subsection{Forms and operators}\label{forms_and_ops}

First of all, we note that 
$$
q_{\eps}(f) = \int_{L(1)}df \wedge \#_{\eps} df = \int_{L(1)}\langle df ,df \rangle_{\eps,W} d\mu_{\Sa}
$$
and that an analogous formula holds for the Sasaki metric. In particular, with Proposition \ref{induced_dual} that implies
$$
q_{\eps}(f) = \int_{L(1)}\left\lbrack\langle df,df\rangle_{\Sa,\eps,W} - \frac{1}{3}\langle df^{\perp},R_W^* df^{\perp}\rangle_{\Sa,\pi(W)} + \eps r_{\eps}(df,df)\right\rbrack d\mu_{\Sa}.
$$
Hence, we obtain:

\begin{Lemma}\label{decomp_ops} For $0<\eps\leq 1$, we have $q_{\eps}(f) = q_{\Sa,\eps}(f) + \Omega(f) + \eps r_{\eps}(f)$, with 
$$
\Omega (f) := - \frac{1}{3}\int_{L(1)}\langle df^{\perp},R_W^* df^{\perp}\rangle_{\Sa,\pi (W)}d\mu_{\Sa}, \, r_{\eps}(f) := \int_{L(1)}r_{\eps}(df,df)d\mu_{\Sa}.
$$
\end{Lemma}

\noindent Let $W\in L(1)$ and $e_{l+1},\ldots e_{m}$ be an orthonormal base of $N_{\pi(W)}L$. We consider vector fields $Z_{\alpha\mu}:\pi^{-1}(x)\to T\pi^{-1}(x)$ given by
\begin{equation}\label{rot_vect_field}
	Z_{\alpha\mu}(W) := J_W (\langle W,e_{\alpha}\rangle\, e_{\mu} - \langle W,e_{\mu}\rangle\, e_{\alpha}) 
\end{equation}
for $l+1\leq \alpha,\mu\leq m$.

\begin{Lemma}\label{ortho} For all sections $W:L\to NL$ and $X:L\to TL$ with horizontal lift $X^h$, we have
$$
\langle J_WW,Z_{\alpha\mu}\rangle_{\Sa} = \langle X^h,Z_{\alpha\mu}\rangle_{\Sa} = \langle J_WW, X^h\rangle_{\Sa} = 0.
$$
\end{Lemma}

\noindent{\bf Proof.} Since $J_WW$ and $Z_{\alpha\mu}$ are vertical vector fields, we have
$$
\langle J_WW,Z_{\alpha\mu}\rangle_{\Sa} = \langle W,K_WZ_{\alpha\mu}\rangle_{NL} = \langle W,e_{\mu}\rangle\langle W,e_{\alpha}\rangle-\langle W,e_{\alpha}\rangle\langle W,e_{\mu}\rangle=0.
$$
The latter equations follow from $X^h\in\mathrm{ker} K_W$ since $X^h$ is horizontal.\hfill $\Box$\\

\noindent Let now $0 < r<1$ and $\partial L (r):= \lbrace W\in NL\,:\, \Vert W\Vert_{NL} = r\rbrace$ the boundary of the $r$-tube around $L$ with Riemannian volume measure $\mu_r$ induced by the Sasaki-metric on $L(1)$. Recall that $\pi:L(1)\to L$ denotes the tube projection. For $x\in L$, we denote the $r$-sphere in $\pi^{-1}(x)$ by $S_{x,r} := \partial L(r)\cap\pi^{-1}(x)$ with induced Riemannian volume $\mu_{x,r}$. The following statement is a direct consequence of Lemma \ref{ortho}.

\begin{Lemma}\label{vflds} {\rm (i)} For all $x\in L$, $0<r<1$, the restriction of the vector fields $Z_{\alpha\mu}$ to $S_{x,r}$ are vector fields on $S_{x,r}$, i.e. $Z_{\alpha\mu}\vert_{S_{x,r}}\in \mathrm{Vect} (S_{x,r})$.
{\rm (ii)} For all $x\in L$, we have $J_WW\vert_{\pi^{-1}(x)}\in\mathrm{Vect}(\pi^{-1}(x))$ 
and for all $0 < r < 1$ and $X\in\mathrm{Vect}(L)$, we have $X^h\vert_{\partial L(r)}\in\mathrm{Vect}(\partial L(r))$. In particular, the restriction of the horizontal bundle $H:=\ker K_W$ to $\partial L(r)$ is a subbundle of $T\partial L(r)$ for all $0<r<1$.
\end{Lemma}

\noindent From this statement, we may conclude that the quadratic forms $q_V, q_H$ and $\Omega$ are decomposable, each one with respect to one of the following three foliations of the tube $L(1)$: 
{\rm{(i)}} $L(1)= \bigcup_{x\in L} \pi^{-1}(x)$, {\rm{(ii)}} $L(1) = L\cup\bigcup_{0< r < 1} \partial L(r)$, and \,finally {\rm{(iii)}} $L(1) = L\cup \bigcup_{0 < r < 1, x\in L} S_{x,r}$. Please note that $L\subset L(1)$ is a zero set, such that we can essentially ignore this part for the discussion of the quadratic forms. \\

\noindent{\bf Remark.} The metric on $\pi^{-1}(x)$ induced by the Sasaki metric is the flat Euclidean metric. The vector fields $Z_{\alpha\mu}$ generate orthogonal transformations and are therefore Killing vector fields on $\pi^{-1}(x)$ and on the spheres $S_{x,r}$, for all $0<r<1$.

\begin{Proposition}\label{decomposedforms} Denote by $\mu_x$, $\mu_{r}$ and $\mu_{x,r}$ the Riemannian volume measures induced by the Sasaki metric on $\pi^{-1}(x)$, $\partial L(r)$ and $S_{x,r}$, respectively. Consider the forms
\begin{enumerate}[label =(\roman*)]
\item[] $q_x (f) := \int_{\pi^{-1}(x)} \langle J^* df, J^* df\rangle_{N^*L} d\mu_x$, 
\item[] $q_r (f) := \int_{\partial L(r)} \langle \kappa df, \kappa df\rangle_{L}d\mu_{r}$,
\item[] $q_{x,r} (f) := \frac{1}{3}\int_{S_{x,r}} \langle df^{\perp},R_W^*df^{\perp}\rangle_{NL}d\mu_{x,r}$,
\end{enumerate}
where $f\in C^{\infty}(L(1))$. Then, the quadratic forms $q_V$, $q_H$ and $\Omega$ are decomposable in the sense that 
\begin{enumerate}[label =(\roman*)]
\item $q_V (f) = \int_L q_x(f\vert_{\pi^{-1}(x)}) d\mu_L$, 
\item $q_H (f) = \int_{(0,1)} q_r(f\vert_{\partial L(r)}) dr$,
\item $\Omega (f) = -\int_{(0,1)\times L} q_{x,r}(f\vert_{S_{x,r}}) dr d\mu_L$.
\end{enumerate}
\end{Proposition}

\noindent{\bf Proof.} Let $x=\pi (W)$ and $e_1, \ldots ,e_l$ an orthonormal base of $T_xL$ and $e_{l+1}, \ldots ,e_m$ an orthonormal base of $N_xL$. As a convention, we denote indices less or equal to $l$ by latin, and larger indices by greek letters. \\
\noindent {\rm (i)} By  
$$
\langle J_W^* df, J_W^* df\rangle_{N^*_xL} = \sum_{\alpha=l+1}^m J_W^*df(e_{\alpha})^2 =  \sum_{\alpha=l+1}^m J_We_{\alpha}(f)^2 
$$
and Lemma \ref{vflds}, {\rm (ii)}, we have $J_We_{\alpha}\in T_W\pi^{-1}(x)$ and therefore, $q_x(f) = q_x(f\vert_{\pi^{-1}(x)})$ for all $x\in L$. {\rm (ii)} By
$$
\langle \kappa df, \kappa df\rangle_{L} = \sum_{j=1}^l \kappa df(e_j)^2 =  \sum_{j=1}^l e_j^h( f)^2
$$
and Lemma \ref{vflds}, {\rm (ii)}, we have $e_j^h\in T\partial L(r)$ and $q_r(f) = q_r(f\vert_{\partial L(r)})$ for all $r\in (0,1)$.
{\rm (iii)} By $W=\sum_{\mu=l+1}^{m}\langle W,e_{\mu}\rangle e_{\mu} = W^{\mu} e_{\mu}\in NL$ (Einstein summation convention), we obtain by the symmetries of the curvature tensor
\begin{eqnarray*}
& & \langle df^{\perp},R_W^*df^{\perp}\rangle_{\Sa,x} \\
&=& \sum_{\alpha = l+1}^{m} df^{\perp}(e_{\alpha}) df^{\perp}(R(W,e_{\alpha})W) \\ &=&\sum_{\alpha,\beta = l+1}^{m} \, J_We_{\alpha}(f) \,\langle R(W,e_{\alpha})W,e_{\beta}\rangle_{NL}\, J_We_{\beta}(f) \\
	&=& \sum_{\alpha, \beta,\mu,\nu = l+1}^{m} \, J_We_{\alpha}(f)W^{\mu}W^{\nu} \,\langle R(e_{\mu},e_{\alpha})e_{\nu},e_{\beta}\rangle_{NL}\, J_We_{\beta}(f) \\
	&=& \frac{1}{4}\sum_{\alpha, \beta,\mu,\nu = l+1}^{m} \,\langle R(e_{\mu},e_{\alpha})e_{\nu},e_{\beta}\rangle_{NL}\,Z_{\alpha\mu}f Z_{\beta\nu}f \\
	&=& \sum_{\alpha < \mu,\beta < \nu} \,\langle R(e_{\mu},e_{\alpha})e_{\nu},e_{\beta}\rangle_{NL}\,Z_{\alpha\mu}f\,Z_{\beta\nu}f
\end{eqnarray*}
with vector fields $Z_{\alpha\mu}$, $Z_{\beta\nu}$ as in (\ref{rot_vect_field}). By Lemma \ref{vflds}, {\rm (i)}, we have $Z_{\alpha\mu}\in TS_{x,r}$, and therefore, $q_{x,r}(f) = q_{x,r}(f\vert_{S_{x,r}})$ for all $r\in (0,1)$, $x\in L$.\hfill$\Box$\\

\noindent{\bf Remark.} (1) The fibre $\pi^{-1}(x)$ with the metric induced from $g_{\Sa}$ is isometric to the flat unit ball. Hence, $q_x$ with domain $\bsob^1(\pi^{-1}(x))$ is the quadratic form of the Dirichlet Laplacian $\Delta_x$ on the flat unit ball. That implies the direct integral decomposition of the vertical operator in (1.1.b). (2) For $r\in (0,1)$, the quadratic form $q_r$ with domain $\mathcal{D}_r$ from Definition \ref{qvqh} {\rm (ii)} is non-negative and closed. By Friedrichs' construction, there is exactly one self-adjoint operator $G_r$ associated to it. The differential expression for $G_r$ is given by $G_r\phi = -\star_r d \star_r \pi^*\kappa d\phi$, where $\star_r$ denotes the Hodge operator associated to the induced metric on $\partial L(r)$.	Therefore, $\Delta_H = \int_{(0,1)}^{\oplus}G_r$, where the operators $G_r$ are self-adjoint and semi-elliptic. \\

\noindent Finally, we collect some facts about the operators associated to the respective quadratic forms which we introduced so far. We will need them in the course of the argument.\\	

\begin{Proposition}\label{operators} The renormalized operator $H_{\eps}^0$ can be written as
$$
H_{\eps}^0 = H_{\Sa,\eps}^0  + P + \eps R_{\eps},
$$
where $R_{\eps}$ is a second order differential expression with smooth coefficients, which are bounded together with all their derivatives uniformly in $1 > \eps > 0$, and 
$$
Pf = \frac{1}{3} \sum_{\alpha < \mu,\beta < \nu} \,\langle R(e_{\mu},e_{\alpha})e_{\nu},e_{\beta}\rangle_{NL}\,Z_{\alpha\mu}\,Z_{\beta\nu}f,
$$
where $e_{\alpha}$, $\alpha = l+1,...,m$ is an arbitrary orthonormal base of $N_{\pi(W)}L$ and $Z_{\alpha\mu}$ is a vector field as in (\ref{rot_vect_field}).
\end{Proposition}

\noindent{\bf Proof.} By partial integration, we obtain, for $f\in\bsob^1\cap\sob^2(L(1),\mu_{\Sa})$, in local coordinates
\begin{eqnarray*}
&& -\frac{1}{3}\sum_{\stackrel{\alpha < \mu}{\beta < \nu}}\int_L \langle R(e_{\mu},e_{\alpha})e_{\nu},e_{\beta}\rangle_{NL}\int_{\pi^{-1}(x)}Z_{\alpha\mu}fZ_{\beta\nu}f\sqrt{\det g_L}dx dw \\
&=& \frac{1}{3}\sum_{\stackrel{\alpha < \mu}{\beta < \nu}}\int_{L(1)}f \langle R(e_{\mu},e_{\alpha})e_{\nu},e_{\beta}\rangle_{NL}Z_{\alpha\mu}Z_{\beta\nu}f d\mu_{\Sa} 
\end{eqnarray*}
by Proposition \ref{decomposedforms} which establishes the statement for $\Omega(f)$. The statement for $R_{\eps}$ follows again by partial integration from the corresponding statement for $r_{\eps}(df,df)$ from Proposition \ref{induced_dual} together with Lemma \ref{decomp_ops}. Hence, 
\begin{eqnarray*}
q_{\eps}^0(f) &=& q_{\Sa,\eps}^0(f) + \Omega(f) + \eps r_{\eps}(f)\\
&=&\int_{L(1)} f\left( H_{\Sa,\eps}^0  + P + \eps R_{\eps}\right)f d\mu_{\Sa},
\end{eqnarray*}
for $f\in\bsob^1\cap\sob^2(L(1),\mu_{\Sa})$, and that implies the statement.
\hfill $\Box$

\begin{Corollary}\label{about_p} {\rm (i)} $\lbrack \Delta_V, P\rbrack = 0$, {\rm (ii)} $PE_0 = 0$.
\end{Corollary}

\noindent{\bf Proof.} {\rm (i)} By Lemma \ref{vflds}, {\rm(i)}, we have
\begin{eqnarray*}
&   &\left.(\lbrack \Delta_V,P\rbrack f)\right\vert_{\pi^{-1}(x)} = \lbrack \Delta_V,P\rbrack (f\vert_{\pi^{-1}(x)}) \\ 
&=& \frac{1}{3} \sum_{\alpha < \mu,\beta < \nu} \,\langle R(e_{\mu},e_{\alpha})e_{\nu},e_{\beta}\rangle_{NL}\lbrack \Delta_V,Z_{\alpha\mu}\,Z_{\beta\nu}\rbrack (f\vert_{\pi^{-1}(x)}) = 0,
\end{eqnarray*}
since 
$$
\lbrack \Delta_V,Z_{\alpha\mu}\,Z_{\beta\nu}\rbrack f\vert_{\pi^{-1}(x)}
= \left(\lbrack\Delta_V,Z_{\alpha\mu}\rbrack Z_{\beta\nu} - Z_{\alpha\mu} \lbrack Z_{\beta\nu}, \Delta_V\rbrack \right)(f\vert_{\pi^{-1}(x)}),
$$
$\Delta_V$ is the Laplacian on every fibre and the $Z_{\alpha\mu}$, $l+1\leq\alpha < \mu\leq m$, are Killing vector fields on the fibres with respect to the Sasaki metric. {\rm (ii)} Let $f\in E_0$. By Proposition \ref{pertprop}, and again by Lemma \ref{vflds},
$$
(Pf)\vert_{\pi^{-1}(x)}= P(f\vert_{\pi^{-1}(x)}) = \overline{f_b} P\phi_0\vert_{\pi^{-1}(x)} = 0,
$$
since $\phi_0$ is invariant under orthogonal transformations of the fibre and the vector fields $Z_{\alpha\mu}$ generate orthogonal transformations.\hfill $\Box$

\subsection{The Proof of Proposition \ref{pertprop}}\label{proof_pertprop}

\noindent Now we are going to prove the different statements of Proposition \ref{pertprop}. \\

\noindent (1) By \ref{sasaki_can_var}.b, $E_0$ is the constant fibre direct integral
$$
E_0 = \int^{\oplus}_L E_{0,x},
$$
where $E_{0,x}$ denotes the (projection onto) the eigenspace corresponding to the lowest eigenvalue $\lambda_0>0$ of the Dirichlet Laplacian on the flat unit ball $B\subset \R^{m-l}$. $E_{0,x}$ is therefore one-dimensional and generated by a normed eigenfunction $\varphi_0$, which is invariant under rotations.  Thus, $\varphi_0(x) = \widetilde{\varphi}(\Vert x\Vert)$. Therefore, it makes sense to define a function $\phi_0\in L^2 (L(1),\mu_{\Sa})$ by $\phi_0(y) :=C\widetilde{\varphi}(d_{\Sa}(y,L))$ where we can choose $C\in\R$ such that $\phi_0$ is non-negative and normalized with respect to the Hilbert space norm.  Thus, a function $f\in E_0$ is determined by a function $f_b : L\to\R$ with 
$$
f_b (x) = \langle f, \phi_0\rangle_{\pi^{-1}(x)}=  \int_{\pi^{-1}(x)}f\phi_0 d\mu_x
$$
where $\mu_x$ denotes the measure on $\pi^{-1}(x)$ which is induced by the volume associated to the Sasaki metric, i.e.
$$
\int_{L(1)}g d\mu_{\Sa} = \int_L \int_{\pi^{-1}(x)}g d\mu_x\,d\mu_L
$$
for all integrable $g$, and $\mu_L$ denotes the Riemannian volume on $L$. Thus, $\overline{f_b}=f_b\circ\pi$ is basic and $f=\overline{f_b}\phi_0$.\\

\noindent (2) Let $f\in E_0\cap\mathcal{D}_H$. Then, by Definition \ref{qvqh}, part (1) above and $\kappa d(\overline{f_b}\phi_0) = \phi_0 \,df_b$, we have
\begin{eqnarray*}
q_H(f) &=& \int_{L(1)} \langle \kappa d(\overline{f_b}\phi_0), \kappa d(\overline{f_b}\phi_0)\rangle_{L}d\mu_{\Sa} = \int_{L(1)} \langle \phi_0 df_b, \phi_0 df_b\rangle_{L}d\mu_{\Sa} \\
&=& \int_{L} \langle df_b, df_b\rangle_{L}\int_{\pi^{-1}(x)} \phi_0^2 d\mu_x \,d\mu_L = \int_{L} \langle df_b, df_b\rangle_{L}\,d\mu_L .
\end{eqnarray*}

\noindent (3) By Lemma \ref{decomp_ops}, we have $l_{\eps}(f) = q_{\eps}(f) - q_{\Sa,\eps}(f) = \Omega(f) + \eps \,r_{\eps}(f).$ \\

\noindent{\rm i.} Since $r_{\eps}(-,-)$ is a bilinear form with smooth and uniformly bounded coefficients, there is a constant $k'_r >0$ such that
$$
\vert r_{\eps,W}(df,df)\vert \leq k'_r\,\langle df,df\rangle_{\Sa,W}
$$
uniformly for all $W\in L(1)$, $0<\eps\leq 1$. That implies by Lemma \ref{est1}, (i)
\begin{eqnarray*}
\vert r_{\eps}(f) \vert &=& \eps\left\vert\int_{L(1)} r_{\eps}(df,df) d\mu_{\Sa}\right\vert \leq \eps k'_r \int_{L(1)} \langle df,df\rangle_{\Sa} d\mu_{\Sa}\\
&\leq& \eps k_r \Vert f\Vert^2_{\bsob^1(L(1),\mu_{\Sa})}.
\end{eqnarray*}

\noindent{\rm ii.} The restriction of the curvature tensor $R$ of $M$ to the submanifold $L$ is a smooth section of the bundle $T^{3,1}M\vert_L$ over $L$. Therefore, the norm of $P_N^*R_W^*P_N^*:N^*_{\pi(x)}L \to N^*_{\pi(x)}L$ is uniformly bounded by some constant $K>0$, i.e. $\Vert P_N^*R_W^*P_N^*\Vert \leq K$ for $W\in NL$. That implies by Cauchy-Schwarz 
$$
\vert\langle df^{\perp},R_W^*df^{\perp}\rangle_{\Sa,\pi (W)} \vert \leq K \,\langle df^{\perp},df^{\perp}\rangle_{\Sa,\pi (W)} = K  \,\langle J_W^*df,J_W^*df\rangle_{N^{*}L}
$$
for $f\in\mathcal{D}_V$. Thus,
\begin{eqnarray*}
\left\vert\Omega (f)\right\vert &=& \frac{1}{3}\left\vert\int_{L(1)} \langle df^{\perp},R_W^*df^{\perp}\rangle_{\Sa,\pi (W)}d\mu_{\Sa}\right\vert \\
&\leq& \frac{K}{3}\int_{L(1)} \langle J_W^*df,J_W^*df\rangle_{N^*L}d\mu_{\Sa} = k_{\Omega} q_V(f).
\end{eqnarray*}

\noindent{\rm iii.} By Proposition \ref{decomposedforms}, (iii), the statement is proved, whenever it is shown for every single fibre. Let thus $x\in L$ and $e_{\alpha}$, $\alpha = l+1,...,m$ be an orthonormal base of $N_xL$. Then
$$
\langle df^{\perp},R_W^*df^{\perp}\rangle_{\Sa,x} = \sum_{\alpha < \mu,\beta < \nu} \,\langle R(e_{\mu},e_{\alpha})e_{\nu},e_{\beta}\rangle_{NL}\,Z_{\alpha\mu}f\,Z_{\beta\nu}f 
$$
with vector fields $Z_{\alpha\mu}$ as in (\ref{rot_vect_field}). The vector field $Z_{\alpha\mu}$ corresponds to the Lie derivative of a one-parameter family of rotations of the plane $e_{\alpha}\wedge e_{\mu}\subset N_xL$. The metric on $\pi^{-1}(x)$ induced by the Sasaki metric is the flat metric. Thus, the vector fields $Z_{\alpha\mu}$ are Killing vector fields of the fibre $\pi^{-1}(x)$. That implies that we have $Z_{\alpha\mu}E_{k,x}\subset E_{k,x}$ for all $l+1\leq\alpha,\mu \leq m$, i.e. the finite-dimensional eigenspaces $E_{k,x}\subset C^{\infty}(\pi^{-1}(x)))$ are invariant under application of the vector fields. Let now $f\in\bsob^1(\pi^{-1}(x))\subset L^2(\pi^{-1}(x))$ and 
$$
f = \sum_{k\geq 0} f_k , \,\, f_k := E_{k,x}f
$$
the orthogonal expansion with smooth $f_k$. Hence 
$$
f= f_0 + (1-E_{0,x})f = f_0 + E_{0,x}^{\perp}f .
$$
By (1), $f_0=a \phi_0$, where $a\in\R$ and $\phi_0$ is invariant with respect to rotations. Hence $Z_{\alpha\mu}f_0 = 0$ for all $l+1\leq\alpha,\mu \leq m$. Thus, $q_{x,r}(f) = q_{x,r}(f_0 + E_{0,x}^{\perp}f) = q_{x,r}(E_{0,x}^{\perp}f)$ and therefore,
\begin{eqnarray*}
\Omega (f) &=& -\int_{(0,1)\times L} q_{x,r}(f\vert_{S_{x,r}}) dr d\mu_L \\
&=& -\int_{(0,1)\times L} q_{x,r}(E_{0,x}^{\perp}f\vert_{S_{x,r}}) dr d\mu_L \\
&=& -\int_{(0,1)\times L} q_{x,r}\left(\left\lbrack E_{0,x}^{\perp}f\right\rbrack\vert_{S_{x,r}}\right) dr d\mu_L = \Omega (E_0^{\perp}f).
\end{eqnarray*}
\hfill $\Box$

\section{Regularity}\label{regularity}

For fixed $\eps > 0$, we consider the set of smooth vectors of $H_{\eps}$ (cf. \cite{ReeSim:75}, X.6, p. 200 ff.) denoted by
$$
\mathcal{C}^{\infty}(H_{\eps}):=\bigcap_{n\geq 1}\mathcal{D}(H_{\eps}^n),
$$
where $\mathcal{D}(H_{\eps}^n)=\lbrace u\in\sob^{2n}(L(1))\,:\,H_{\eps}^{n-1}u\vert_{\partial L(1)} = ... = u\vert_{\partial L(1)}=0\rbrace$. The sets $\mathcal{C}^{\infty}(H_{\eps})$ and $\mathcal{C}^{\infty}(H_{\eps}^0)$ coincide.\\

\noindent{\bf Remark.} Let $f\in L^2(L(1),\mu_{\Sa})$. By \cite{Lun:95}, Proposition 2.1.1 {\rm(i)}, we have 
$$
e^{-\frac{t}{2}H_{\eps}^0}f \in \mathcal{C}^{\infty}(H_{\eps})
$$
for all $t>0$.\\

\noindent We are now going to consider different norms on the set of smooth vectors to finally prove that the semigroups generated by $H_{\eps}^0$ actually converge smoothly in the sense of Theorem \ref{Main}.

\subsection{Boundary conditions}

By examining the boundary conditions and by considering smooth vectors as solutions of another boundary problem for which we have elliptic a priori estimates, we are going to construct a family of norms which are equivalent to $2n$-Hilbert Sobolev norms on the set of smooth vectors. Let $\Delta_{\Sa}$ be the Laplacian on $L(1)$ associated to the Sasaki metric and $\Vert - \Vert_{2n}$ the $2n$-Sobolev norm on $\sob^{2n}(L(1),\mu_{\Sa})$.

\begin{Lemma}\label{boundary_conditions} Let $u\in \CC^{\infty}( H_{\eps}^0)$. For all $n\geq 1$, there are differential operators $T_{n}(\eps)$, $1\geq \eps > 0$, defined in a neighbourhood of $L(1)$, such that
	\begin{enumerate}[label = (\roman*)]
		\item $\mathrm{ord} \, T_n(\eps) \leq 2n$, 
		\item all coefficients are smooth and bounded together with their derivatives uniformly in $1\geq \eps > 0$,
		\item we have $\Delta_V^n u\vert_{\partial L(1)} = \eps^3\, T_n(\eps)u\vert_{\partial L(1)}$.
	\end{enumerate}
\end{Lemma}

\noindent{\bf Proof.} {\rm a.} The case $n=1$ is provided by 
	\begin{eqnarray*}
		0 &=& H_{\eps}^0u\vert_{\partial L(1)} = \eps^2\,  H_{\eps}^0u\vert_{\partial L(1)} \\
		&=&(\Delta_V -\lambda_0 + \eps^2(\Delta_H + P)  + \eps^3 R_{\eps})u\vert_{\partial L(1)} \\
		&=&(\Delta_V + \eps^3 R_{\eps})u\vert_{\partial L(1)} ,
	\end{eqnarray*}
	because $\Delta_Hu\vert_{\partial L(1)}=Pu\vert_{\partial L(1)}=0$ since the restrictions of these operators to $\partial L(1)$ yield proper differential operators on $\partial L(1)$ where $u$ is constant (and equal to zero). By Proposition \ref{operators}, the operator $T_1(\eps) := - R_{\eps}$ satisfies {\rm (i), (ii)}.
	{\rm b.} Assume now that the induction hypothesis is valid for $1,...,n-1$. Then, keeping in mind that $\Delta_V$ commutes with $\Delta_H$ and $P$ on smooth functions,
	\begin{eqnarray*}
		0 &=& \left(H_{\eps}^0\right)^nu\vert_{\partial L(1)} = \left(\eps^2\,  H_{\eps}^0\right)^n u\vert_{\partial L(1)} \\
		&=&(\Delta_V -\lambda_0 + \eps^2(\Delta_H + P)  + \eps^3 R_{\eps})^nu\vert_{\partial L(1)} \\
		&=&((\Delta_V -\lambda_0)^n + n\eps^2(\Delta_V -\lambda_0)^{n-1}(\Delta_H + P)  + \eps^3 G_{\eps})u\vert_{\partial L(1)} \\
		&=&(\Delta_V^n + \sum_{r=0}^{n-1}\left(\begin{array}{c}n \\ r\end{array}\right) (-\lambda_0)^{n-r}\Delta_V^r )u\vert_{\partial L(1)}  + \eps^3 G_{\eps}u\vert_{\partial L(1)} \\ 
		&&+ n\eps^2\sum_{r=0}^{n-1}\left(\begin{array}{c}n-1 \\ r\end{array}\right) (-\lambda_0)^{n-r-1}(\Delta_H + P)\Delta_V^ru\vert_{\partial L(1)}, 
	\end{eqnarray*} 
where $G_{\eps}$ is a differential operator of order less or equal to $2n$ with uniformly bounded smooth coefficients.
	From the induction hypothesis together with Corollary \ref{about_p}, all remaining terms can be summarized to an operator $T_n(\eps)$ satisfying conditions {\rm (i)}, {\rm (ii)}, by
	\begin{eqnarray*}
		0 &=&\left(\Delta_V^n + \eps^3\sum_{r=1}^{n-1}\left(\begin{array}{c}n \\ r\end{array}\right) (-\lambda_0)^{n-r}T_r(\eps)\right)u \vert_{\partial L(1)}  + \eps^3 G_{\eps}u\vert_{\partial L(1)} \\ 
		&&+ n\eps^5\sum_{r=1}^{n-1}\left(\begin{array}{c}n-1 \\ r\end{array}\right) (-\lambda_0)^{n-r-1}(\Delta_H + P)T_r(\eps)u\vert_{\partial L(1)} \\
		&=:&\Delta_V^nu\vert_{\partial L(1)}  - \eps^3T_n(\eps)u\vert_{\partial L(1)} .  
	\end{eqnarray*} 
	That implies the statement.\hfill$\Box$\\

\noindent From this result, we immediately conclude the corresponding result for $\Delta_{\Sa}$.

\begin{Corollary}\label{b_c} Let $u\in \CC^{\infty}( H_{\eps}^0)$. For all $n\geq 1$, there are differential expressions $S_{n}(\eps)$, $1\geq \eps > 0$, defined in a neighbourhood of $L(1)$ such that
	\begin{enumerate}[label=(\roman{*})]
		\item $\mathrm{ord} \, S_n(\eps) \leq 2n$, 
		\item all coefficients are smooth and bounded with all their derivatives uniformly in $1\geq \eps > 0$,
		\item we have $\Delta_{\Sa}^n u\vert_{\partial L(1)} = \eps^3\, S_n(\eps)u\vert_{\partial L(1)}$.
	\end{enumerate}
\end{Corollary}

\noindent{\bf Proof.} We have, again by \ref{sasaki_can_var}.c and Lemma \ref{boundary_conditions} that
	\begin{eqnarray*}
		\Delta_{\Sa}^nu\vert_{\partial L(1)} &=& \sum_{r=0}^{n}\left(\begin{array}{c}n \\ r\end{array}\right) \Delta_H^{n-r}\Delta_V^ru \vert_{\partial L(1)}  \\
		&=& \Delta_H^n u\vert_{\partial L(1)} + \sum_{r=1}^{n}\left(\begin{array}{c}n \\ r\end{array}\right) \Delta_H^{n-r}(\Delta_V^r-\eps^3 T_r(\eps))u \vert_{\partial L(1)}  \\
		&&+ \eps^3\sum_{r=1}^{n}\left(\begin{array}{c}n \\ r\end{array}\right) \Delta_H^{n-r} T_r(\eps)u \vert_{\partial L(1)}.
	\end{eqnarray*}
	Since $u\vert_{\partial L(1)}=0$ implies $\Delta_H^n u\vert_{\partial L(1)}=0$, we obtain
	$$
	\Delta_{\Sa}^nu \vert_{\partial L(1)} = \eps^3\sum_{r=1}^{n}\left(\begin{array}{c}n \\ r\end{array}\right) \Delta_H^{n-r} T_r(\eps)u \vert_{\partial L(1)} =: \eps^3\,S_n(\eps)u \vert_{\partial L(1)},
	$$
	and by Lemma \ref{boundary_conditions}, the operator $S_n(\eps)$ is indeed defined in a neighbourhood of $L(1)$. \hfill$\Box$\\

\noindent Therefore, we may consider functions $u\in\CC^{\infty}( H_{\eps}^0)$ as solutions of the boundary problem
\begin{equation}\label{sasabound}
\left\lbrace\begin{array}{l}
\Delta_{\Sa}^nu = f(u)\,\,\,\, \mathrm{on}\,\, L(1), \\  \\
\Delta_{\Sa}^{n-1}u\vert_{\partial L(1)} =  \eps^3\, S_{n-1}(\eps)u\vert_{\partial L(1)},\\
\,\,\,\,\,\,\,\,\,\,\,\,\,\vdots\\
\Delta_{\Sa}\,\,\,\,u\vert_{\partial L(1)} =  \eps^3\,S_1(\eps)u\vert_{\partial L(1)},\\
u\vert_{\partial L(1)} = 0.
\end{array}\right. .
\end{equation}
Here, $f(u) := \Delta_{\Sa}^n u\in C^{\infty}(L(1))$. This boundary problem satisfies the Shapiro - Lopatinskij conditions (cf. \cite{Agr:97}, Sect. 1.3, p. 8 ff) and is therefore regular elliptic. Hence, we obtain the following elliptic a priori estimate (cf. \cite{Agr:97}, Thm. 2.2.1, p. 16) as a first alternative representation of the $2n$-Sobolev norm on the space of smooth vectors.

\begin{Proposition}\label{esti2} For every $n\geq 1$ there is some $C_n > 0$ and $\eps_n > 0$ such that for all $0<\eps < \eps_n$, we have
	\begin{equation}\label{apriori2}
	\Vert u\Vert_{2n} \leq C_n \left( \Vert u\Vert  + \Vert\Delta_{\Sa}^nu\Vert \right) 
	\end{equation}
	for all $u\in \CC^{\infty}( H_{\eps}^0)$.
\end{Proposition}

\noindent{\bf Proof.} The elliptic a priori estimate (\cite{Agr:97}, Thm. 2.2.1, p. 16), applied to (\ref{sasabound}) above, reads
	\begin{eqnarray*}
		\Vert u\Vert_{2n} &\leq& C'_n \left( \Vert u\Vert  + \Vert\Delta_{\Sa}^nu\Vert  + \sum_{r=1}^{n-1}\Vert \Delta_{\Sa}^ru\Vert_{\mathrm{H}^{2n-2r-1/2}(\partial L(1))}\right) \\
		&=&C'_n \left( \Vert u\Vert + \Vert\Delta_{\Sa}^nu\Vert + \eps^3\sum_{r=1}^{n-1}\Vert S_r(\eps)u\Vert_{\mathrm{H}^{2n-2r-1/2}(\partial L(1))}\right) \\
		&=&C'_n\left( \Vert u\Vert + \Vert\Delta_{\Sa}^nu\Vert + \eps^3c_n\Vert u\Vert_{2n}\right) 
	\end{eqnarray*}
	since, by construction, the operators $S_r(\eps)$ from Lemma \ref{b_c} are globally defined differential expressions on $L(1)$. For $\eps > 0$ small enough, we may absorb the last term on the right hand side into the left hand side and obtain the statement after redefining the constant. \hfill$\Box$

\subsection{A scale of norms on smooth vectors}

\noindent 
We are going to establish a family of norms which are equivalent to $2n$-Hilbert Sobolev norms on the set of smooth vectors. Let $\Delta_{\Sa}$ be the Laplacian on $L(1)$ associated to the Sasaki metric and $\Vert - \Vert_{2n}$ the $2n$-Sobolev norm on $\sob^{2n}(L(1),\mu_{\Sa})$.

\subsubsection{Uniform regularity: The case $\mathbf{n=1}$.}

We first treat the case $n=1$ and prove a Kato-type inequality for the solution of the boundary value problem which yields an estimate of the $2$-Sobolev norm of a function $u\in\DD( H_{\eps}^0)$.

\begin{Proposition}\label{ind_anf} Let $u\in \CC^{\infty}( H_{\eps}^0)$. Then, there is an $\eps_0 > 0$ and a constant $B_1 > 0$ such that
$$
\Vert u\Vert_2 \leq B_1 \left(\Vert u \Vert + \Vert  H_{\eps}^0 u\Vert \right) 
$$
uniformly for all $0<\eps < \eps_0$.
\end{Proposition}  

\noindent{\bf Proof.} By the spectral theorem, we have
\begin{eqnarray*}
&&\Vert H_{\Sa,\eps}^0u\Vert^2 \\
&=& \frac{1}{\eps^4} \Vert(\Delta_V-\lambda_0)u\Vert^2 + \Vert \Delta_Hu\Vert^2 + \frac{2}{\eps^2}\langle \Delta_Hu, (\Delta_V - \lambda_0)u\rangle\\
&=& \Vert \Delta_Hu\Vert^2 + \sum_{k\geq 1}\left\lbrack\frac{1}{\eps^2}\left(\frac{\lambda_k - \lambda_0}{\eps}\right)^2\Vert E_ku\Vert^2 + \frac{2}{\eps}\frac{\lambda_k - \lambda_0}{\eps}\langle\Delta_Hu,E_ku\rangle\right\rbrack
\end{eqnarray*}
and, letting $\eps \leq 1-\lambda_0/\lambda_1$, we obtain $\eps^{-1}(\lambda_k - \lambda_0)\geq \lambda_k$ for all $k\geq 1$. Thus, using the shorthand $u^{\perp} = (1-E_0)u$ and the fact that by $\lbrack E_k,\Delta_H\rbrack = 0$ (see \ref{sasaki_can_var}.c), the last summand
$$
\langle \Delta_Hu,E_ku\rangle = \langle \Delta_Hu,E_k^2u\rangle = \langle \Delta_HE_ku,E_ku\rangle \geq 0
$$
is non - negative, we obtain
\begin{equation}\label{in1}
\Vert H_{\Sa,\eps}^0u\Vert^2 \geq \frac{1}{\eps^2} \Vert \Delta_Vu^{\perp}\Vert^2 + \Vert \Delta_Hu\Vert^2 + \frac{2}{\eps}\langle \Delta_Hu^{\perp},\Delta_Vu^{\perp}\rangle .
\end{equation}
Since all three summands are non - negative, this particularly implies 
\begin{equation}\label{inequality}
\Vert H_{\Sa,\eps}^0u\Vert \geq \left\lbrace\begin{array}{c}\frac{1}{\eps} \Vert \Delta_Vu^{\perp}\Vert \\ \Vert \Delta_Hu\Vert \\ \sqrt{\frac{2}{\eps}}\langle \Delta_Hu^{\perp},\Delta_Vu^{\perp}\rangle^{1/2}\end{array}\right. .
\end{equation}
Now, by Propositon \ref{decomposedforms}, the quadratic form and $\Omega$ and therefore also the associated operator $P$ is decomposable with respect to the direct integral decomposition of $L^2(L(1),\mu_{\Sa})$ by the Hilbert spaces on the fibres. By Corollary \ref{about_p}, {\rm (ii)}, we have $Pu = Pu^{\perp}$.  Since $P$ is a second order differential operator with bounded coefficients that implies 
\begin{eqnarray*}
\Vert Pu \Vert^2 &=& \int_L \Vert P (u^{\perp}\vert_{\pi^{-1}(x)})\Vert_x^2 \, d\mu_L \\
&\leq& \int_L c^2 \left(\Vert u^{\perp}\vert_{\pi^{-1}(x)}\Vert_x^2  + \Vert \Delta_V (u^{\perp}\vert_{\pi^{-1}(x)})\Vert_x^2 \right)\, d\mu_L \\
&\leq& c^2\left( \Vert u^{\perp} \Vert +\Vert \Delta_V u^{\perp}\Vert\right)^2.
\end{eqnarray*}
By Proposition \ref{operators}, $H_{\eps}^0 = H_{\Sa,\eps}^0 + P  + \eps R_{\eps}$. First of all, we have by (\ref{inequality}) and $H_{\Sa,\eps}^0 u^{\perp} = H_{\Sa,\eps}^0 u$,
$$
\Vert Pu \Vert\leq c \left( \Vert u\Vert + \Vert \Delta_V u^{\perp}\Vert\right)\leq c \left( \Vert u\Vert + \eps\Vert H_{\Sa,\eps}^0 u\Vert\right) .
$$
On the other hand, the remainder is of second order and satisfies $\Vert R_{\eps}u\Vert\leq d\,\Vert u\Vert_2$. Hence, for some $b>0$, we have
\begin{eqnarray*}
& & \Vert H_{\eps}^0u\Vert + b \Vert u\Vert \\
&=& \Vert(H_{\Sa,\eps}^0 + P  + \eps 
R_{\eps})u\Vert + b \Vert u\Vert \\
&\geq& \Vert H_{\Sa,\eps}^0u\Vert - \Vert P u\Vert -\eps\Vert R_{\eps}u\Vert + b \Vert u\Vert \\
&\geq& \Vert H_{\Sa,\eps}^0u\Vert - c\left(\Vert u\Vert + \eps\Vert H_{\Sa,\eps}^0u\Vert\right)-\eps\,d\Vert u\Vert_2 + b \Vert u\Vert .
\end{eqnarray*}
By $\eps < 1-\lambda_0/\lambda_1<1$, inequality (\ref{in1}) together with $\lbrack \Delta_V,\Delta_H\rbrack =0$ implies 
$$
\Vert H_{\Sa,\eps}^0u\Vert^2 \geq\Vert \Delta_Vu^{\perp}\Vert^2 + \Vert \Delta_Hu\Vert^2 + 2\langle \Delta_Hu,\Delta_Vu^{\perp}\rangle = \Vert\Delta_Vu^{\perp} +  \Delta_Hu\Vert^2,
$$
and therefore,
$$
\Vert H_{\Sa,\eps}^0u\Vert \geq \Vert\Delta_{\Sa}u + \Delta_V(u^{\perp}-u)\Vert\geq \Vert\Delta_{\Sa}u \Vert - \lambda_0\Vert E_0u\Vert\geq\Vert\Delta_{\Sa}u \Vert - \lambda_0\Vert u\Vert.
$$
Thus, we have
\begin{eqnarray*}
\Vert H_{\eps}^0u\Vert + b \Vert u\Vert 
&\geq& (1-c\eps) \Vert H_{\Sa,\eps}^0u\Vert + (b - c)\Vert u\Vert -\eps\,d\Vert u\Vert_2 \\
&\geq& (1-c\eps) \Vert\Delta_{\Sa}u\Vert + (b - c - \lambda_0(1-c\eps))\Vert u\Vert -\eps\,d\Vert u\Vert_2 . 
\end{eqnarray*}
By Proposition \ref{esti2}, we have
$$
\Vert u\Vert_2 \leq C_1\left( \Vert u\Vert + \Vert\Delta_{\Sa}u\Vert\right)
$$
and hence,
$$
\Vert H_{\eps}^0u\Vert + b \Vert u\Vert 
\geq (1-(c+dC_1)\eps) \Vert\Delta_{\Sa}u\Vert + (b - c-\lambda_0+(c\lambda_0 - dC_1)\eps)\Vert u\Vert    . 
$$
Now we choose $\eps < \eps_0$ small enough to obtain $k_1:= 1-(c+dC_1)\eps > 0$ and, at the same time, $b > 0$ large enough to obtain $k_2 := b - c-\lambda_0+(c\lambda_0 - dC_1)\eps > 0$. Thus, again by Proposition \ref{esti2},
$$
\Vert H_{\eps}^0u\Vert + b \Vert u\Vert 
\geq k_1\Vert\Delta_{\Sa}u\Vert + k_2\Vert u\Vert \geq \frac{\min \lbrace k_1,k_2\rbrace}{C_1}\,\Vert u \Vert_2 
$$
and we finally obtain the statement by letting
$B_1 := \frac{C_1\,\max\lbrace 1,b\rbrace}{\min \lbrace k_1,k_2\rbrace}$.
\hfill $\Box$

\subsubsection{Uniform regularity: The estimate for $\mathbf{n > 1}$.} 

We will now derive another alternative representation of the $2n$-Sobolev norm of an element $u\in\CC^{\infty}( H_{\eps}^0)$. We are going to use the following consequence of the Calderon-Lions interpolation theorem (\cite{ReeSim:75}, Theorem IX.20, p. 37), which we state without proof.

\begin{Proposition}\label{interpol} Let $1\leq r\leq 2n-1$ be an integer. Then, for every $1\geq \sigma > 0$ there are constants $C(n),C(\sigma,n)>0$ such that
	$$
	\Vert u\Vert_r \leq C(n)\left(\sigma \Vert u\Vert_{2n} + C(\sigma,n)\Vert u\Vert\right).
	$$
\end{Proposition} 

\noindent Now we come to the result just announced. The assertion is proved by induction.

\begin{Proposition}\label{ind_schlu} Let $u = u(\eps) \in\CC^{\infty}( H_{\eps}^0)$. Then, for all $n\geq 1$ there are numbers $B_n > 0$, $\eps_n > 0$, such that
$$
\Vert u\Vert_{2n} \leq B_n \left(\Vert u\Vert +\Vert H_{\eps}^0u\Vert_{2n-2}\right)
$$
uniformly for $0<\eps < \eps_n$.
\end{Proposition}

\noindent{\bf Proof.} The case $n=1$ was already treated in Proposition \ref{ind_anf} and will be used in the course of the argument. By Proposition \ref{esti2}, we have 
$$
\Vert u\Vert_{2n} \leq C_n \left( \Vert\Delta_{\Sa}^nu\Vert + \Vert u\Vert\right) \leq C_n'\left( \Vert\Delta_{\Sa}^{n-1}u\Vert_2 + \Vert u\Vert\right),
$$
since $\Delta_{\Sa} : \mathrm{H}^{m+2}(L(1),\mu_{\Sa})\to\mathrm{H}^{m}(L(1),\mu_{\Sa})$ is continuous. Recall that the operators $S_n(\eps)$ above are actually defined on an open neighbourhood of $\overline{L(1)}$. Hence,
$$
\Vert u\Vert_{2n}  \leq C_n'\left( \Vert\Delta_{\Sa}^{n-1}u-\eps^3S_{n-1}(\eps)u\Vert_2 +\eps^3 \Vert S_{n-1}(\eps)u\Vert_2+ \Vert u\Vert\right)
$$
and since by construction $\Delta_{\Sa}^{n-1}u-\eps^3S_{n-1}(\eps)u\in\bsob^1\cap\mathrm{H}^{2}(L(1),\mu_{\Sa})$, we may use Proposition \ref{ind_anf} and obtain with $C_n'' := C_n'B_1$
\begin{eqnarray*}
\Vert u\Vert_{2n}  &\leq& C_n'' \left(\Vert\Delta_{\Sa}^{n-1}u-\eps^3S_{n-1}(\eps)u\Vert + \Vert H_{\eps}^0(\Delta_{\Sa}^{n-1}u-\eps^3S_{n-1}(\eps)u)\Vert  \right) \\
&& + C_n'\left(\eps^3 \Vert S_{n-1}(\eps)u\Vert_2+ \Vert u\Vert\right) \\
&\leq& A_1 \Vert u\Vert_{2n-2} + A_2 \Vert H_{\eps}^0\Delta_{\Sa}^{n-1}u\Vert +A_3\eps \Vert \eps^2  H_{\eps}^0S_{n-1}(\eps)u\Vert \\
&&+ A_4 \eps^3 \Vert u\Vert_{2n} + A_5 \Vert u \Vert
\end{eqnarray*}
with constants $A_1,...,A_5 > 0$. Since $\eps^2 H_{\eps}^0$ is a smooth differential expression with all coefficients bounded uniformly in $\eps > 0$ together with their derivatives, we have
$$
A_3\eps \Vert \eps^2  H_{\eps}^0S_{n-1}(\eps)u\Vert \leq A_3' \eps\Vert u \Vert_{2n}.
$$
By Proposition \ref{interpol}, taking $\sigma > 0$ small enough, we may now absorb the term with Sobolev index $2n-2$ into the left hand side. Furthermore, by taking $\eps > 0$ small enough, we may absorb the term with Sobolev index $2n$ as well. Thus, after changing the constants accordingly, we obtain an $a_n>0$ with
$$
\Vert u\Vert_{2n} \leq a_n \left(\Vert u\Vert + \Vert H_{\eps}^0\Delta_{\Sa}^{n-1}u\Vert \right) .
$$
Now, $\lbrack H_{\Sa,\eps}^0,\Delta_{\Sa}\rbrack u = 0$ and hence
$$
\lbrack  H_{\eps}^0,\Delta_{\Sa}^{n-1}\rbrack u = \lbrack P  + \eps R_{\eps},\Delta_{\Sa}^{n-1}\rbrack u=: Q_n(\eps)u,
$$
and $Q_n(\eps)$ is a differential expression of order at most $2n-1$ with coefficients that are uniformly bounded together with all their derivatives independent of $\eps > 0$. Thus
\begin{eqnarray*}
\Vert u\Vert_{2n} &\leq& a_n \left\lbrace\Vert u\Vert + \Vert\Delta_{\Sa}^{n-1} H_{\eps}^0u\Vert  + \Vert Q_n(\eps)u\Vert\right\rbrace \\
&\leq& a_n \left\lbrace\Vert u\Vert + \Vert\Delta_{\Sa}^{n-1} H_{\eps}^0u\Vert  + b_n\Vert u\Vert_{2n-1}\right\rbrace
\end{eqnarray*}
and absorbing again the term involving the $(2n-1)$-Sobolev norm into the left hand side by interpolation using Proposition \ref{interpol}, we obtain the statement, since $\Delta_{\Sa}^{n-1}: \mathrm{H}^{2n-2}(L(1),\mu_{\Sa})\to L^2(L(1),\mu_{\Sa})$ is continuous. \hfill $\Box$\\

\noindent From this result, we derive the last estimate of the Sobolev norm in terms of the operator $ H_{\eps}^0$. Note that the estimate holds uniformly for a family $u(\eps)$ but that we omit the argument in the statement of the estimate. 

%\begin{Corollary} Let $u = u(\eps) \in\CC^{\infty}( H_{\eps}^0)$. Then for all $n\geq 1$ there are numbers $D_n > 0$, $\eps_n > 0$, such that
%\begin{equation}\label{abschatz}
%\Vert u\Vert_{2n} \leq D_n \sum_{k=0}^n \Vert  \left( H_{\eps}^0\right)^ku\Vert
%\end{equation}
%uniformly for $\eps < \eps_n$.
%\end{Corollary}
%
%\noindent{\bf Proof.} {\rm a.} The case $n=1$ is provided by Proposition \ref{ind_anf} with $D_1=B_1$. {\rm b.} For $n > 1$, we have by Proposition \ref{ind_schlu}  $\Vert u\Vert_{2n} \leq B_n \left(\Vert H_{\eps}^0u\Vert_{2n-2} + \Vert u\Vert\right)$ . Note now that $u\in \CC^{\infty}( H_{\eps}^0)$ also implies $ H_{\eps}^0u\in \CC^{\infty}( H_{\eps}^0)$ and that if we assume that the statement is valid for $1,...,n-1$, we obtain
%$$
%\Vert  H_{\eps}^0u\Vert_{2n-2} \leq D_{n-1} \sum_{k=0}^{n-1} \Vert  \left( H_{\eps}^0\right)^{k+1}u\Vert
%$$
%and hence by Proposition \ref{ind_schlu}
%\begin{eqnarray*}
%\Vert u\Vert_{2n} &\leq& B_n\left(D_{n-1} \sum_{k=0}^{n-1} \Vert   \left( H_{\eps}^0\right)^{k+1}u\Vert + \Vert u\Vert\right) \\
%&\leq& B_n\left( D_{n-1} \sum_{k=1}^{n} \Vert  \left( H_{\eps}^0\right)^{k}u\Vert + \Vert u\Vert\right)\\
%&\leq&  D_{n} \sum_{k=0}^{n} \Vert  \left( H_{\eps}^0\right)^{k}u\Vert 
%\end{eqnarray*}
%with suitably chosen $D_n > 0$. That implies the statement.\hfill $\Box$

\begin{Corollary}\label{convenient} Let $u = u(\eps) \in\CC^{\infty}( H_{\eps}^0)$ and $\alpha \geq 2\max\lbrace \lambda_0,1\rbrace$. Then, for all $n\geq 1$, there are numbers $D_n > 0$ and $\eps_n > 0$, such that
\begin{equation}\label{abschatz}
\Vert u\Vert_{2n} \leq D_n \left(\Vert u\Vert + \Vert  \left( H_{\eps}^0+\alpha\right)^nu\Vert\right),
\end{equation}
uniformly for $\eps < \eps_n$.
\end{Corollary}

\noindent{\bf Proof.} By Corollary \ref{spc}, $\mathrm{spec} (H_{\eps}^0 + \alpha)\subset \lbrack 2\max\lbrace \lambda_0,1\rbrace-\lambda_0,\infty)\subset \lbrack 1,\infty)$ uniformly for $0<\eps<\eps_0$. That implies $\lambda_s(\eps) + \alpha \geq \vert \lambda_s(\eps)\vert$ for the eigenvalues $\lambda_s(\eps)$ of the operators $H_{\eps}^0$ and therefore,
$$
(\lambda_s(\eps) + \alpha )^k \geq (\lambda_s(\eps) + \alpha)^l \geq \lambda_s(\eps)^l,
$$
for all $1\leq l\leq k$, $\eps < \eps_0$ and $s\geq 0$. By the spectral theorem that implies
$$
\Vert (H_{\eps}^0 )^l u\Vert\leq \Vert \left( H_{\eps}^0+\alpha\right)^lu\Vert \leq \Vert\left( H_{\eps}^0+\alpha\right)^ku\Vert,
$$
for all $1\leq l\leq k$, $\eps < \eps_0$ and $s\geq 0$. We now proceed by induction: {\rm a.} The case $n=1$ is provided by Proposition \ref{ind_anf} and
$$
\Vert u\Vert_{2} \leq B_1 \left(\Vert u\Vert + \Vert H_{\eps}^0u\Vert\right)\leq B_1 \left(\Vert u\Vert + \Vert  \left( H_{\eps}^0+\alpha\right)u\Vert\right),
$$
hence $D_1=B_1$. {\rm b.} Assume now
$$
\Vert  u\Vert_{2n-2}  \leq D_{n-1} \left(\Vert  u\Vert + \Vert  \left( H_{\eps}^0+\alpha\right)^{n-1} u\Vert\right),
$$
for $n - 1 \geq 1$. By Proposition \ref{ind_schlu}  $\Vert u\Vert_{2n} \leq B_n \left(\Vert u\Vert + \Vert H_{\eps}^0u\Vert_{2n-2}\right)$ uniformly for $0<\eps\leq\eps_n$. Note now that $u\in \CC^{\infty}( H_{\eps}^0)$ also implies $ H_{\eps}^0u\in \CC^{\infty}( H_{\eps}^0)$ and therefore
$$
\Vert  H_{\eps}^0u\Vert_{2n-2} \leq D_{n-1} \left(\Vert H_{\eps}^0u\Vert + \Vert  \left( H_{\eps}^0+\alpha\right)^{n-1}H_{\eps}^0u\Vert\right).
$$
By the spectral theorem
\begin{eqnarray*}
\Vert \left( H_{\eps}^0+\alpha\right)^{n-1}H_{\eps}^0u\Vert^2 &=& \sum_{s\geq 0} \left( \lambda_s(\eps)+\alpha\right)^{2n-2}\lambda_s(\eps)^2\vert\langle u_s(\eps),u\rangle\vert^2 \\
&\leq& \sum_{s\geq 0} \left( \lambda_s(\eps)+\alpha\right)^{2n}\vert\langle u_s(\eps),u\rangle\vert^2 \\
&\leq& \Vert \left( H_{\eps}^0+\alpha\right)^{n}u\Vert^2,
\end{eqnarray*}
where for $s\geq 0$, $\lambda_s(\eps)$ and $u_s(\eps)$ denote eigenvalues and eigenfunctions of $ H_{\eps}^0$, respectively. Hence, by Proposition \ref{ind_schlu}
\begin{eqnarray*}
\Vert u\Vert_{2n} &\leq& B_n\left(\Vert u\Vert + D_{n-1} \left(\Vert H_{\eps}^0u\Vert +  \Vert   \left( H_{\eps}^0+\alpha\right)^{n-1}H_{\eps}^0u\Vert \right)\right) \\
&\leq& B_n\left(\Vert u\Vert + D_{n-1} \left(2  \Vert   \left( H_{\eps}^0+\alpha\right)^{n}u\Vert \right)\right)\\
&\leq&  D_{n} \left(\Vert  u\Vert + \Vert  \left( H_{\eps}^0+\alpha\right)^n u\Vert\right)
\end{eqnarray*}
with $D_n := B_n \,\max\lbrace 1,2D_{n-1}\rbrace$. That implies the statement.\hfill $\Box$

\subsection{Convergence of the semigroups.}

Now we can finally use the representation (\ref{abschatz}) of the $2n$-Sobolev norm to apply the spectral theorem. Let now $u(\eps):\lbrack 0,1\rbrack\to L^2(L(1),\mu_{Sa})$ be a strongly continuous family. Let furthermore $\eps_n,\alpha >0$ as in Corollary \ref{convenient}. For $\eps,t>0$, we consider now
$$
u(\eps,t) := \exp\left( -\frac{t}{2}\left(H_{\eps}^0 +\alpha\right) \right) u(\eps),
$$
and for $\overline{\eps},\overline{t}>0$, we denote by $\UU(\overline{\eps},\overline{t})$ the set
$$
\UU (\overline{\eps},\overline{t}) := \lbrace u(\eps,t) \,:\, 0<\eps\leq \overline{\eps}, t\geq \overline{t}\rbrace .
$$

\begin{Proposition}\label{combo} {\rm (i)} $\UU (\overline{\eps},\overline{t}) \subset\CC^{\infty}( H_{\eps}^0)$. {\rm (ii)} For all $n\geq 0$, the subset $\UU (\eps_n,\overline{t})\subset \mathrm{H}^{2n}(L(1),\mu_{\Sa})$ is uniformly bounded.  {\rm (iii)} $\UU (\eps_n,\overline{t})$ is equi - continuous as a subset of 
$$
C(\lbrack \overline{t},\infty),\mathrm{H}^{2n}(L(1),\mu_{\Sa})),
$$ 
i.e. for all $\gamma >0$ there is some $\delta >0$ such that $\vert t-t'\vert < \delta$ implies $\Vert u(\eps,t) - u(\eps,t')\Vert_{2n}\leq \gamma$ uniformly in $0 < \eps \leq \eps_n$.
\end{Proposition}

\noindent{\bf Proof.} {\rm (i)} This is a basic property of analytic semigroups, cf. \cite{Lun:95}, Proposition 2.1.1 {\rm (i)}, p. 35. {\rm (ii)} By the spectral theorem, we have
$$
\Vert \left( H_{\eps}^0+\alpha\right)^{n}u(\eps,t)\Vert^2 = \sum_{s\geq 0} (\lambda_s(\eps)+\alpha)^{2n}e^{-t(\lambda_s(\eps)+\alpha)}\left\vert\langle u_s(\eps),u(\eps) \rangle\right\vert^2,
$$
where for $s\geq 0$, $\lambda_s(\eps)$ and $u_s(\eps)$ denote again the eigenvalues and eigenfunctions of $ H_{\eps}^0$. Note that by Corollary \ref{spc}, $\lambda_s(\eps)+\alpha \geq 1>0$ for all $0 < \eps\leq \eps_0$ and all $n\geq 0$. Therefore, by $x^{2n}e^{-tx} \leq \left(\frac{2n}{et}\right)^{2n}$ for $x\geq 0$, we obtain since $\Vert u(\eps)\Vert$ is uniformly bounded by the continuity assumption
$$
\Vert \left( H_{\eps}^0 + \alpha\right)^{n}u(\eps,t)\Vert^2 \leq \left(\frac{2n}{et}\right)^{2n}\Vert u(\eps) \Vert^2\leq \left(\frac{2n}{e\overline{t}}\right)^{2n}\Vert u(\eps) \Vert^2 =:M_n^2.
$$
Boundedness now follows from inequality (\ref{abschatz}).
{\rm (iii)} By
$$
\left\vert e^{-bt'} - e^{-bt''} \right\vert\leq  b e^{-b \overline{t}}\,\vert t'-t''\vert ,
$$
for all $t',t''\in\lbrack \overline{t},\infty)$ and $b >0$, we obtain with the same estimate as above
\begin{eqnarray*}
&&\Vert  \left( H_{\eps}^0+\alpha \right)^{n}(u(\eps,t') - u(\eps,t''))\Vert^2 \\
&=& \sum_{s\geq 0} (\lambda_s(\eps)+\alpha)^{2n}\left\vert e^{-\frac{t'(\lambda_s(\eps) +\alpha) }{2}}-e^{-\frac{t''(\lambda_s(\eps)+\alpha )}{2}}\right\vert^2\left\vert\langle u_s(\eps),u(\eps) \rangle\right\vert^2\\
&\leq& \frac{1}{4}\vert t' - t''\vert^2\, \sum_{s\geq 0} e^{-\overline{t}(\lambda_s(\eps)+\alpha)}(\lambda_s(\eps)+\alpha)^{2n+2}\left\vert\langle u_s(\eps),u(\eps) \rangle\right\vert^2\\
&\leq& \frac{1}{4}\vert t' - t''\vert^2\,  \left(\frac{2n}{e\overline{t}}\right)^{2n}\Vert u(\eps) \Vert^2  \\
&\leq& \frac{M_n^2}{4} \vert t' - t''\vert^2.
\end{eqnarray*}
By (\ref{abschatz}), that implies equi - continuity.\hfill $\Box$

\subsection{Proof of Theorem \ref{Main}} {\rm a.} Let $\overline{t} > 0$. By Proposition \ref{combo}, {\rm (ii)}, the set 
$$
\UU (t) := \lbrace u(\eps,t) \,:\,\eps\leq\eps_{n+1} \rbrace \subset \mathrm{H}^{2n+2}(L(1),\mu_{\Sa})
$$
is uniformly bounded for all $t\geq \overline{t}$ fixed. By the Sobolev embedding theorem, it is therefore relatively compact as a subset of $\sob^{2n}(L(1),\mu_{\Sa})$. By Proposition \ref{combo}, {\rm (iii)}, the set $\UU(\eps_{n+1},\overline{t})\subset C(\lbrack \overline{t},\infty),\sob^{2n}(L(1),\mu_{\Sa})$ is equi - continuous for all $\overline{t} > 0$. Hence, for each compact interval $I\subset (0,\infty)$, the subset
$$
\UU(\eps_{n+1},I):= \lbrace u(\eps,-) : I\to \sob^{2n}(L(1),\mu_{\Sa}), \eps\leq\eps_{n+1}\rbrace
$$
is relatively compact in $C(I,\sob^{2n}(L(1),\mu_{\Sa}))$ by the Arzela - Ascoli theorem.

\noindent {\rm b.} Let now $(\eta_{m})_{m\geq 1}$ with $\eps_0\geq\eta_{m}\geq 0$ be a sequence tending to zero as $m$ tends to infinity and let $(u_m = u(\eta_m,-))_{m\geq 1}\subset \UU(\eps_{n+1}, I)$ be an arbitrary sequence. By {\rm a.}, it contains a subsequence $u_m'$ that converges to some limit 
$$
u_{\infty}'\in C(I,\sob^{2n}(L(1),\mu_{\Sa})).
$$
By
$$
\lim_{m \to\infty}\sup_{t\in I}\Vert u_m' - u_{\infty}'\Vert \leq \lim_{m \to\infty}\sup_{t\geq t_0}\Vert u_m' - u_{\infty}'\Vert_{2n} = 0
$$
and the $L^2$-result from Proposition \ref{semi1} that implies
$$
u_{\infty}'(t) = E_0 e^{-\frac{t}{2}\left(\Delta_L+\alpha\right)}E_0u(0) .
$$
Hence, every subsequence contains a convergent subsequence with the same limit. That implies 
$$
\lim_{\eps \to 0} e^{-\frac{t}{2}\left( H_{\eps}^0 + \alpha\right)}u(\eps) = E_0 e^{-\frac{t}{2}\left(\Delta_L + \alpha\right) }E_0u(0)
$$
in $\sob^{2n}(L(1),\mu_{\Sa})$. Multiplication by $e^{\frac{t\alpha}{2}}$ yields the statement.\hfill $\Box$

\end{document}